\documentclass[11pt,draft]{article}

\usepackage[sort]{cite}       
\usepackage{amssymb}          
\usepackage{latexsym}         
\usepackage{amsmath}          
\usepackage{amscd}            
\usepackage[latin1]{inputenc} 
\usepackage{eucal}            
\usepackage{bbm}              
\usepackage{theorem}          

\title{Convergence of the Wick Star Product}

\author{\textbf{Svea Beiser}\thanks{E-mail:
    Svea.Beiser@physik.uni-freiburg.de}\addtocounter{footnote}{2},
  \textbf{Hartmann Römer}\thanks{E-mail:
    Hartmann.Roemer@physik.uni-freiburg.de}\addtocounter{footnote}{0},
  \textbf{Stefan Waldmann}\thanks{E-mail:
    Stefan.Waldmann@physik.uni-freiburg.de}
  \\[0.1cm]
  Fakult{\"a}t f{\"u}r Mathematik und Physik\\
  Albert-Ludwigs-Universit{\"a}t Freiburg\\
  Physikalisches Institut\\
  Hermann Herder Stra{\ss}e 3\\
  D 79104 Freiburg\\
  Germany}

\date{June 2005\\[0.5cm] FR-THEP 2005/05}

\renewcommand{\mathbb}[1]{\mathbbm{#1}} 

\newcounter{comment}

\newcommand{\cc}[1]      {\overline{{#1}}}

\newcommand{\image}      {\operatorname{\mathrm{im}}}    
    
\newcommand{\Ad}         {\operatorname{\mathrm{Ad}}}

\newcommand{\SP}[1]      {\left\langle{#1}\right\rangle} 
\newcommand{\Unit}       {\mathbb{1}}

\newcommand{\norm}[1]    {\left\|{#1}\right\|}            
 
\newcommand{\I}          {\mathrm{i}}
\newcommand{\E}          {\mathrm{e}}
\newcommand{\D}          {\operatorname{\mathrm{d}}} 

\newcommand{\Exp}        {\operatorname{\mathrm{Exp}}}
\newcommand{\IM}         {\operatorname{\mathrm{Im}}}

\newcommand{\Binom}[2]   {\genfrac{(}{)}{0pt}{1}{#1}{#2}}

\newcommand{\AntiHol}    {\cc{\mathcal{O}}}
\newcommand{\HolAntiHol} {\mathcal{O}\!\!\times\!\!\cc{\mathcal{O}}}

\newcommand{\starwick}  {\mathbin{\star_{\scriptscriptstyle\mathrm{Wick}}}}
\newcommand{\starwickh} {\mathbin{\star_{\scriptscriptstyle\mathrm{Wick}}^\hbar}}
\newcommand{\starwicka} {\mathbin{\star_{\scriptscriptstyle\mathrm{Wick}}^\alpha}}

\newcommand{\hph}[2][m,\ell,R,S]    {h^{p,\hbar}_{#1}\left(#2\right)}
\newcommand{\normph}[2][m,\ell,R,S] {\norm{{#2}}^{p, \hbar}_{#1}}
\newcommand{\normNullh}[2][m,\ell,R,S] {\norm{{#2}}^{0, \hbar}_{#1}}

\newcommand{\tAph}    {\widetilde{\mathcal{A}}_{p,\hbar}}
\newcommand{\Aph}     {\mathcal{A}_{p,\hbar}}
\newcommand{\Ah}      {\mathcal{A}_{\hbar}}
\newcommand{\HBF}     {\mathfrak{H}_{\mathrm{\scriptscriptstyle BF}}}

\newtheorem{lemma} {Lemma} [section]
\newtheorem{proposition} [lemma] {Proposition}

\newtheorem{theorem} [lemma] {Theorem}
\newtheorem{corollary} [lemma] {Corollary}
\newtheorem{definition}[lemma] {Definition}

\theorembodyfont{\rmfamily}
\newtheorem{example}[lemma]{Example}
\newtheorem{remark}[lemma]{Remark}

\newenvironment{proof}[1][{}]{
  \par\noindent
  \textsc{Proof{#1}:}
}
{
  \hspace*{\fill} $\blacksquare$\newline
}

\numberwithin{equation}{section}

\begin{document}

\maketitle

\begin{abstract}
    We construct a Fr\'echet space as a subspace of
    $C^\omega(\mathbb{C}^n)$ where the Wick star product converges and
    is continuous.  The resulting Fr\'echet algebra $\Ah$ is studied
    in detail including a $^*$-representation of $\Ah$ in the
    Bargmann-Fock space and a discussion of star exponentials and
    coherent states.
\end{abstract}

%
%


%
%

\section{Introduction}
\label{sec:Intro}

Deformation quantization usually comes in two flavours: formal and
strict deformations:

In formal deformation quantization as introduced by
\cite{bayen.et.al:1978a}, see also \cite{dito.sternheimer:2002a,
  gutt:2000a} for recent reviews, one considers the Poisson algebra of
smooth complex-valued functions $C^\infty(M)$ on a Poisson manifold as
observable algebra in classical mechanics. A formal star product
$\star$ is a $\mathbb{C}[[\lambda]]$-bilinear associative product for
$C^\infty(M)[[\lambda]]$ such that in zeroth order $f \star g$ is the
pointwise product and in first order of $\lambda$ the
$\star$-commutator gives $\I$ times the Poisson bracket. Usually one
requires the higher orders to be given by bidifferential operators.
The algebra $C^\infty(M)[[\lambda]]$ then serves as a model for the
quantum mechanical observables corresponding to the classical system
described by $M$. In particular, the formal parameter $\lambda$
corresponds to Planck's constant $\hbar$ and should be replaced by
$\hbar$ whenever one can establish convergence of the formal series.

On the other hand, in strict deformation quantization as introduced in
\cite{rieffel:1989a}, see also \cite{landsman:1998a}, one works on
the framework of $C^*$-algebras where the dependence of the deformed
product $\star_\hbar$ on the deformation parameter $\hbar$ is now
required to be continuous. This is made precise using the notion of
continuous fields of $C^*$-algebras.

While in the first approach one has very strong existence
\cite{kontsevich:2003a, dewilde.lecomte:1983b, fedosov:1986a,
  fedosov:1989a, fedosov:1996a,omori.maeda.yoshioka:1991a} and
classification results \cite{kontsevich:2003a, nest.tsygan:1995a,
  nest.tsygan:1995b, bertelson.cahen.gutt:1997a, weinstein.xu:1998a,
  deligne:1995a, gutt.rawnsley:1999a}, the
formal character of Planck's constant is, of course, physically not
acceptable. Here the second approach is much more appealing as it
directly uses the analytical framework suitable for quantum mechanics.
On the other hand, however, a general construction and reasonable
classification of strict quantizations seems still to be missing.

Many examples like the global symbol calculus on cotangent bundles
\cite{bordemann.neumaier.pflaum.waldmann:2003a,
  bordemann.neumaier.waldmann:1999a,
  bordemann.neumaier.waldmann:1998a}, Berezin-Toeplitz quantization on
Kähler manifolds \cite{bordemann.meinrenken.schlichenmaier:1991a,
  karabegov.schlichenmaier:2001a, cahen.gutt.rawnsley:1995a,
  cahen.gutt.rawnsley:1994a, cahen.gutt.rawnsley:1993a,
  cahen.gutt.rawnsley:1990a} as well as \cite{bieliavsky:2002a}
suggest that the formal star products should be seen as asymptotic
expansions for $\hbar \to 0$ of their convergent counterparts in
strict quantization. On the other hand, many formal star products
allow for large subalgebras, where the formal series actually
converge, whence in some sense the asymptotics can be used again to
recover the strict result, a heuristic statement for which a general
theorem unfortunately is still missing. The above examples also
suggest that there is a framework in between formal and
$C^*$-algebraic, namely one can try to construct deformations of
$C^\infty(M)$ (or suitable subalgebras of it) in the framework of
\emph{Fr\'echet} or more generally \emph{locally convex algebras}.
Early results in this direction have been obtained in
\cite{hansen:1984a, kammerer:1986a, maillard:1986a}, see also
\cite{omori.maeda.miyazaki.yoshioka:2000a}. Moreover, a general set-up
of smooth deformations has been established and exemplified in
\cite{dubois-violette.kriegl.maeda.michor:2001a}, in
\cite{pflaum.schottenloher:1998a} holomorphic deformations were
studied.

The example we are going to discuss will provide an entire holomorphic
deformation of a Fr\'echet subalgebra of $C^\omega(M)$.  More
specifically, we consider the most simple phase space $M =
\mathbb{C}^n$ with its canonical Poisson structure $\{z^k,
\cc{z}^\ell\} = \frac{2}{\I} \delta^{k\cc{\ell}}$ and the \emph{formal
  Wick star product}
\begin{equation}
    \label{eq:WickStarProductNull}
    f \starwick g 
    = \sum_{r=0}^\infty \frac{(2\lambda)^r}{r!}
    \sum_{i_1, \ldots, i_r}
    \frac{\partial^r f}
    {\partial z^{i_1} \cdots \partial z^{i_r}}
    \frac{\partial^r g}
    {\partial \cc{z}^{i_1} \cdots \partial \cc{z}^{i_r}},
\end{equation}
where $z^1, \ldots, z^n$ are the canonical, global, holomorphic
coordinates on $\mathbb{C}^n$. Our convergence scheme to construct the
`convergent' subalgebra $\Ah$ is then based on the crucial observation
that the Wick star product enjoys a very strong positivity property
\cite{bursztyn.waldmann:2004a:pre, bursztyn.waldmann:2000a,
  bordemann.waldmann:1998a}: every $\delta$-functional is a positive
$\mathbb{C}[[\lambda]]$-linear functional in the sense of formal power
series. After choosing a point $p \in \mathbb{C}^n$ and $\hbar > 0$
this will allow us to construct recursively a system of seminorms for
which the Wick star product is continuous, thereby defining our
algebra $\Aph$. Moreover, we shall construct, via the GNS construction
corresponding to the positive functional $\delta_p$, a faithful
$^*$-representation of $\Aph$ on a dense subspace of the Bargmann-Fock
space giving an interpretation of the $\delta$-functionals as coherent
states with respect to the Heisenberg group $\mathsf{H}_n$ acting on
$\mathbb{C}^n$. We treat this example in quite some detail as we
believe that it may serve as a good starting point for geometric
generalizations to Wick star products on Kähler manifolds
\cite{karabegov:1996a, bordemann.waldmann:1997a} suitable for a
bottom-up approach to \cite{cahen.gutt.rawnsley:1995a,
  cahen.gutt.rawnsley:1994a, cahen.gutt.rawnsley:1993a,
  cahen.gutt.rawnsley:1990a,
  bordemann.meinrenken.schlichenmaier:1991a,
  karabegov.schlichenmaier:2001a}. Moreover, in a future project we
shall discuss the field-theoretic generalization for infinitely many
degrees of freedom.

The paper is organized as follows: In Section~\ref{sec:preliminary} we
briefly recall some basic properties of $\starwick$, the formal GNS
construction for $\delta$-functionals and the Bargmann-Fock space.
Section~\ref{sec:Construction} is devoted to the construction of the
seminorms, depending on $p \in \mathbb{C}^n$ and $\hbar > 0$. This
gives the space $\Aph$ which is shown to be a subspace of
$C^\omega(\mathbb{C}^n)$ with a Fr\'echet topology. In
Section~\ref{sec:ContinuityWick} we show that $\Aph$ is a subalgebra
of $C^\omega(\mathbb{C})$ such that the pointwise product as well as
the Poisson bracket are continuous, i.e.  $\Aph$ becomes a
Fr\'echet-Poisson algebra. Moreover, we show that the formula
\eqref{eq:WickStarProductNull} for $\starwick$ actually converges on
$\Aph$ in the Fr\'echet topology resulting in a continuous product.
This way, $\Aph$ becomes a holomorphic deformation. In
Section~\ref{sec:TranslationRescale} we discuss the dependence on the
a priori chosen point $p$ and on the value $\hbar > 0$ of Planck's
constant. It turns out that the translation group acts on $\Aph$ by
inner $^*$-automorphisms whence $\Aph = \Ah$ does not depend on the
choice of $p$. As a side remark we show that the star exponential, see
\cite{bayen.et.al:1978a}, of linear functions converges in the
topology of $\Ah$. Moreover, for all values $\hbar > 0$ the algebras
$\Ah$ are isomorphic in a canonical way. Finally, in
Section~\ref{sec:GNSCoherentStates} we show how the GNS construction
yields a $^*$-representation of $\Ah$ in the sense of
\cite{schmuedgen:1990a} in the Bargmann-Fock space. The action of the
Heisenberg group by inner $^*$-automorphisms gives easily the coherent
states.

\medskip

\noindent
\textbf{Acknowledgement:} We would like to thank Pierre Bieliavsky,
Nico Giulini, Simone Gutt, Nikolai Neumaier, Konrad Schmüdgen, Martin
Schlichenmaier and Rainer Verch for valuable discussions on this
example. Moreover, S.~W. thanks the ULB for hospitality while part of
this work was being done.

%
%

\section{Preliminary results}
\label{sec:preliminary}

In this section we shall collect some well-known results on the Wick
star product which we shall use in the sequel, see e.g.
\cite{bordemann.waldmann:1998a, bordemann.waldmann:1997a}.

On the classical phase space $\mathbb{R}^{2n} \cong \mathbb{C}^n$ with
standard symplectic form $\omega = \frac{\I}{2} \D z^k \wedge \D
\cc{z}^k$ one defines the \emph{formal Wick star product} by
\begin{equation}
    \label{eq:WickStarProduct}
    f \starwick g =
    \sum_{r=0}^\infty \frac{(2\lambda)^r}{r!}
    \sum_{i_1, \ldots, i_r} 
    \frac{\partial^r f}
    {\partial z^{i_1} \cdots \partial z^{i_r}}
    \frac{\partial^r g}
    {\partial \cc{z}^{i_1} \cdots \partial \cc{z}^{i_r}},
\end{equation}
where $f, g \in C^\infty(\mathbb{C}^n)[[\lambda]]$, the formal
parameter $\lambda$ corresponds to Planck's constant $\hbar$ without
any further prefactors and $z$, $\cc{z}$ denote the usual global
holomorphic/anti-holomorphic coordinates on $\mathbb{C}^n$. Then
$\starwick$ is known to be an associative star product quantizing the
canonical Poisson bracket corresponding to $\omega$. It has the
separation of variable property in the sense of Karabegov
\cite{karabegov:2000a, karabegov:1996a} and is in fact the name-giving
example of a star product of Wick type in the sense of
\cite{bordemann.waldmann:1997a}.  Moreover, $\starwick$ is Hermitian
\begin{equation}
    \label{eq:WickHermitian}
    \cc{f \starwick g} = \cc{g} \starwick \cc{f},
\end{equation}
where according to our interpretation of $\lambda$ the formal
parameter is defined to be real $\cc{\lambda} = \lambda$. We shall
also frequently make use of multiindex notation: Let $R = (r_1,
\ldots, r_n) \in \mathbb{N}^n$ be a multiindex, then one defines $|R|
= r_1 + \cdots + r_n$, $R! = r_1! \cdots r_n!$ as well as $z^R =
(z^1)^{r_1} \cdots (z^n)^{r_n}$ etc. Moreover, we define $R \le L$ if
$r_i \le \ell_i$ for all $i=1, \ldots, n$. The Wick star product can
equivalently be written as
\begin{equation}
    \label{eq:WickMultiindex}
    f \starwick g = 
    \sum_{R=0}^\infty \frac{(2\lambda)^{|R|}}{R!}
    \frac{\partial^{|R|} f}{\partial z^R}
    \frac{\partial^{|R|} g}{\partial \cc{z}^R}.
\end{equation}

The Wick star product enjoys a very strong positivity property which
e.g.  the Weyl-Moyal star product does not share: If $\delta_p:
C^\infty(\mathbb{C}^n)[[\lambda]] \longrightarrow
\mathbb{C}[[\lambda]]$ denotes the evaluation functional at $p \in
\mathbb{C}^n$ then we have
\begin{equation}
    \label{eq:DeltapWick}    
    \delta_p(\cc{f} \starwick f)
    =
    \sum_{R=0}^\infty \frac{(2\lambda)^{|R|}}{R!}
    \left|\frac{\partial^{|R|} f}{\partial \cc{z}^R}(p)\right|^2
    \ge 0,
\end{equation}
where the positivity is understood in the sense of the canonical ring
ordering of $\mathbb{R}[[\lambda]]$, see
\cite{bordemann.waldmann:1998a, waldmann:2005b} for a detailed
discussion on the physical relevance of this notion of positivity.
This very strong positivity property is not true for general Hermitian
star products: instead one has to add `quantum corrections' to a given
classically positive functional (here $\delta_p$) in order to obtain a
positive functional with respect to the star product. Since positive
functionals play the role of \emph{states} the above simple
observation implies that for the Wick star product any classical state
defines a quantum state \emph{without any quantum correction}, see
\cite{waldmann:2005b} for a detailed discussion on the general
situation. In fact, the Wick star product is used in an essential way
for proving that for an arbitrary Hermitian star product one can
always construct quantum corrections for a classical state, see
\cite{bursztyn.waldmann:2004a:pre, bursztyn.waldmann:2000a}.

Since $\delta_p$ is a positive functional for $\starwick$ we have a
corresponding GNS representation of the Wick star product algebra
$C^\infty(\mathbb{C}^n)[[\lambda]]$. In fact, this example was one of
the first examples of GNS constructions in deformation quantization.
Since we need the construction in the following, we briefly review the
results from \cite{bordemann.waldmann:1998a}. The Gel'fand ideal
$\mathcal{J}_p$ of $\delta_p$ is given by
\begin{equation}
    \label{eq:GelfandIdeal}
    \mathcal{J}_p = 
    \left\{
        f \; \Big| \;
        \delta_p( \cc{f} \starwick f) = 0
    \right\}
    =
    \left\{
        f \; \left| \;
            \frac{\partial^{|R|}f}{\partial \cc{z}^R}(p) = 0
            \; \textrm{for all} \; R
        \right.
    \right\}.
\end{equation}
The GNS pre Hilbert space $\mathcal{H}_p =
C^\infty(\mathbb{C}^n)[[\lambda]] \big/ \mathcal{J}_p$ can canonically
be identified with
\begin{equation}
    \label{eq:FormalBFSpace}
    \mathcal{H}_{\mathrm{\scriptscriptstyle{BF}}}
    =
    \mathbb{C}[[\cc{y}^1, \ldots, \cc{y}^n]][[\lambda]]
\end{equation}
with the $\mathbb{C}[[\lambda]]$-valued positive definite inner
product
\begin{equation}
    \label{eq:FormalBFInnerProduct}
    \SP{\phi, \psi}
    =
    \sum_{R=0}^\infty \frac{(2\lambda)^{|R|}}{R!}
    \; \cc{\frac{\partial^{|R|} \phi}{\partial \cc{y}^R}(0)}
    \; \frac{\partial^{|R|} \psi}{\partial \cc{y}^R}(0),
\end{equation}
where the identification of a class $\psi_f \in \mathcal{H}_p$ is
given by its formal anti-holomorphic Taylor expansion, i.e.
\begin{equation}
    \label{eq:IdentGNSHpFormalBF}
    \mathcal{H}_p \ni \psi_f 
    \; \mapsto \;
    \sum_{R=0}^\infty \frac{1}{R!} 
    \; \frac{\partial^{|R|}f}{\partial \cc{z}^R}(p) 
    \; \cc{y}^R,
\end{equation}
where $\psi_f$ denotes the equivalence class of the function $f$ in
$\mathcal{H}_p$. Then the GNS representation on $\mathcal{H}_p$
defined by $\pi_p(f)\psi_g = \psi_{f \starwick g}$ is translated into
\begin{equation}
    \label{eq:FormalBFRepAtp}
    \varrho_p(f) = \sum_{R, S=0}^\infty
    \frac{(2\lambda)^{|R|}}{R!S!}
    \; \frac{\partial^{|R| + |S|}f}{\partial z^R \partial \cc{z}^S}(p)
    \; \cc{y}^S 
    \; \frac{\partial^{|R|}}{\partial \cc{y}^R},
\end{equation}
via the unitary map \eqref{eq:IdentGNSHpFormalBF}. In particular, for
$p=0$ we see that this gives the formal analog of the usual
\emph{Bargmann-Fock space} and the \emph{Bargmann-Fock
  representation}: indeed, recall that the Bargmann-Fock space is the
Hilbert space
\begin{equation}
    \label{eq:BargmannFockSpace}
    \HBF = 
    \left\{
        f \in \AntiHol(\mathbb{C}^n) 
        \; \left| \;
            \frac{1}{(2\pi\hbar)^n}
            \int |f(\cc{z})|^2 \E^{-\frac{\cc{z}z}{2\hbar}} \D z \D \cc{z}
            < \infty
        \right.
    \right\}
\end{equation}
of \emph{anti-holomorphic} functions which are square-integrable with
respect to the Gaussian measure, see \cite{bargmann:1961a}.  Then it
is well-known that $\HBF$ is actually a closed subspace of the
$L^2$-space for this measure and hence a Hilbert space itself.
Moreover, the $L^2$-inner product can be evaluated by the same formula
\eqref{eq:FormalBFInnerProduct} if one replaces $\lambda$ by $\hbar$,
where the series now converges absolutely. A Hilbert basis for $\HBF$
is given by the monomials
\begin{equation}
    \label{eq:BFHilbertBasis}
    \mathsf{e}_{R} (\cc{z}) 
    = \frac{1}{\sqrt{(2\hbar)^{|R|} R!}}  
    \; \cc{z}^R.
\end{equation}
From \eqref{eq:FormalBFInnerProduct} we see that the Hilbert space
$\HBF$ can be interpreted as the space of those anti-holomorphic
functions whose Taylor coefficients at $0$ form a sequence in a
(weighted) $\ell^2$-space.

%
%

\section{Construction of the Fr\'echet space $\Aph$}
\label{sec:Construction}

The motivation for our convergence scheme is rather simple: we fix
$\hbar > 0$ and we fix a point $p \in \mathbb{C}^n$. Then we are
looking for a subalgebra of $C^\infty(\mathbb{C}^n)[[\lambda]]$ such
that $\delta_p(f)$ converges for $\lambda = \hbar$. Though this looks
rather innocent at the beginning, we obtain a hierarchy of conditions:
For $f, g$ in our \emph{algebra} we want $f \starwick g$ to be in the
algebra as well whence $\delta_p(f \starwick g)$ has to converge for
$\lambda = \hbar$ as well. The idea is now to estimate the convergence
of $\delta_p(f \starwick g)\big|_{\lambda = \hbar}$ by
$\delta_p(\cc{f} \starwick f)\big|_{\lambda = \hbar}$ and
$\delta_p(\cc{g} \starwick g)\big|_{\lambda = \hbar}$ using the
Cauchy-Schwarz inequality for the \emph{positive} $\delta$-functional.
An iteration of this procedure will lead to countably many
\emph{unary} conditions for a function $f$ to belong to our algebra,
which we shall interpret as seminorms determining the algebra. This
approach can hence be used as a heuristic motivation for the following
definition of the seminorms.
\begin{definition}
    \label{definition:hmlRS}
    Let $R, S \in \mathbb{N}^n$ be multiindices, $m \in \mathbb{N}$
    and $\ell = 0, \ldots, 2^m - 1$. Then we define recursively for $f
    \in C^\infty(\mathbb{C}^n)$
    \begin{equation}
        \label{eq:hNullNull}
        \begin{split}
            \hph[0,0,R,S]{f}
            &=
            \left.
                (2\hbar)^{|R| + |S|} 
                \left(
                    \frac{\partial^{|R|+|S|}\cc{f}}
                    {\partial z^S \partial \cc{z}^R} 
                    \starwick
                    \frac{\partial^{|R|+|S|}f}
                    {\partial \cc{z}^S \partial z^R}
                \right)(p) 
            \right|_{\lambda = \hbar} \\
            &=
            \sum_{N=0}^\infty \frac{(2\hbar)^{|R| + |S| + |N|}}{N!}
            \left|
                \frac{\partial^{|R| + |S| + |N|} f}
                {\partial z^R \partial \cc{z}^{N+S}}(p)
            \right|^2
        \end{split}
    \end{equation}
    and
    \begin{equation}
        \label{eq:hmlRS}
        \hph{f}
        =
        \begin{cases}
            \displaystyle
            \sum_{N=0}^\infty \frac{1}{N!} 
            \left(
                \sum_{I=0}^R\sum_{J=0}^{N+S}
                \Binom{R}{I}\Binom{N+S}{J}
                \hph[m-1,\ell/2,I,J]{f}
            \right)^2
            & \ell \quad \textrm{even} \\
            \displaystyle
            \sum_{N=0}^\infty \frac{1}{N!} 
            \left(
                \sum_{I=0}^R\sum_{J=0}^{N+S}
                \Binom{R}{I}\Binom{N+S}{J}
                \hph[m-1,(\ell-1)/2,J,I]{f}
            \right)^2
            & \ell \quad \textrm{odd.}
        \end{cases}
    \end{equation}
\end{definition}
The precise form of the recursive definition will become either clear
by following the above heuristic argument in detail or from the proof
of Proposition~\ref{proposition:fWickPolynomialCont}: the binomial
coefficients arise from the Leibniz rule.

Thanks to the positivity \eqref{eq:DeltapWick} of the
$\delta$-functional it is clear that $\hph[0, 0, R, S]{f}$ either
converges absolutely or diverges absolutely to $+ \infty$, as it is a
series consisting of non-negative terms only. By induction, the same
is true for all other $\hph{f}$. Hence we have
\begin{equation}
    \label{eq:hphPositive}
    \hph{f} \in [0, + \infty],
\end{equation}
where `convergence' is always absolute and does not depend on the
order of summation.
\begin{definition}
    \label{definition:Seminormph}
    For $f \in C^\infty(\mathbb{C}^n)$ we define
    \begin{equation}
        \label{eq:SemiNormDef}
        \normph{f} = \sqrt[2^{m+1}]{\hph{f}} 
        \quad \in [0, + \infty].
    \end{equation}
    Moreover, we define
    \begin{equation}
        \label{eq:AtildeDef}
        \tAph = 
        \left\{
            f \in C^\infty(\mathbb{C}^n) 
            \; \Big| \;
            \normph{f} < \infty
            \quad
            \textrm{for all}
            \quad
            m, \ell, R, S
        \right\}.
    \end{equation}
\end{definition}

The first step is now to show that $\tAph$ is a vector space and that
the $\normph{\cdot}$ are seminorms on $\tAph$. This will be a
consequence of the following proposition:
\begin{proposition}
    \label{proposition:SeminormsAreSeminorms}
    The maps $\normph{\cdot}: C^\infty(\mathbb{C}^n) \longrightarrow
    [0, +\infty]$ enjoy the following properties:
    \begin{enumerate}
    \item $\normph{\alpha f} = |\alpha| \normph{f}$ for $\alpha \in
        \mathbb{C}$.
    \item $\normph{f + g} \le \normph{f} + \normph{g}$.
    \item $\normph[m-1,\ell,R,S]{f} \le \normph[m,2\ell,R,S]{f}$ and
        $\normph[m-1,\ell,R,S]{f} \le \normph[m,2\ell+1,R,S]{f}$.
    \item $\normph[m,\ell,0,S]{f} \le \sqrt[2^{m+2}]{S!}
        \normph[m+1,2\ell,0,0]{f}$.
    \item $\normph[m,\ell,R,0]{f} \le
        \sqrt[2^{m+2}]{R!}\normph[m+1,2\ell+1,0,0]{f}$.
    \item $\normph{f} \le \sqrt[2^{m+3}]{R!}\sqrt[2^{m+2}]{S!}
        \normph[m+2,4\ell+1,0,0]{f}$.
    \end{enumerate}
\end{proposition}
\begin{proof}
    The first part is clear by a simple induction. For the second part
    the case $m=0$ follows directly from Minkowski's inequality.  Then
    $m > 0$ is shown inductively by using again Minkowski's
    inequality, for both cases of odd and even $\ell$. The remaining
    inequalities between the $\normph{f}$ for different values of the
    parameters are simply obtained by omitting all but one term for a
    specific $N$ in the defining summations of \eqref{eq:hmlRS}. For
    example, in the third part one considers $N = 0$ and $I = R$, $J =
    S$ only, while for the fourth and fifth one uses $N = S = J$ and
    $N = R = I$.
\end{proof}

Thus the maps $\normph{\cdot}$, restricted to $\tAph$, give indeed
seminorms. Moreover, the labels $R$, $S$ play only a minor role thanks
to the estimates in the last part of the proposition. This motivates
the following definitions. For $f \in C^\infty(\mathbb{C}^n)$ we
define
\begin{equation}
    \label{eq:normmellDef}
    \normph[m, \ell]{f} = \normph[m, \ell, 0, 0]{f}
\end{equation}
\begin{equation}
    \label{eq:normmDef}
    \normph[m]{f} 
    = \max_{0 \le \ell \le 2^m-1}
    \left\{\normph[m,\ell]{f}\right\}.
\end{equation}
Then we have the following simple corollary:
\begin{corollary}
    \label{corollary:TildeAphIsLocallyConvex}
    The set $\tAph \subseteq C^\infty(\mathbb{C}^n)$ is a subvector
    space and the $\normph{\cdot}$ are seminorms on $\tAph$. Moreover,
    the seminorms $\normph[m,\ell]{\cdot}$ as well as the seminorms
    $\normph[m]{\cdot}$ determine the same locally convex topology on
    $\tAph$.
\end{corollary}

In the following, we shall equip $\tAph$ always with this locally
convex topology induced by the seminorms $\normph{\cdot}$. However,
this topology has one unpleasant feature: it is \emph{non-Hausdorff}
as a function $f$ whose Taylor expansion at $p$ vanishes identically
has clearly $\normph{f} = 0$ for all parameters $m, \ell, R, S$. On
the other hand, as one sees already from the seminorm
$\normph[1]{\cdot}$ the functions with \emph{vanishing} $\infty$-jet
$\mathsf{j}^\infty_pf$ at $p$ are the only functions with this
property. Thus we identify them to be zero in order to have a
Hausdorff topology:
\begin{definition}
    \label{definition:Aph}
    We define
    \begin{equation}
        \label{eq:AphDef}
        \Aph = 
        \tAph \big/
        \left\{
            f \in C^\infty(\mathbb{C}^n) 
            \; \big| \; 
            \mathsf{j}^\infty_pf = 0
        \right\},
    \end{equation}
    and equip $\Aph$ with the induced locally convex topology
    determined by the seminorms $\normph{\cdot}$ (or equivalently, by
    the seminorms $\normph[m]{\cdot}$).
\end{definition}

Clearly, $\Aph$ is now a Hausdorff locally convex topological vector
space. The following theorem shows that the abstract quotient can be
viewed as a certain subspace of the real-analytic functions on
$\mathbb{C}^n$:
\begin{theorem}
    \label{theorem:RealAnalytic}
    Let $\hbar > 0$ and $p \in \mathbb{C}^n$.
    \begin{enumerate}
    \item $\Aph$ is a Hausdorff locally convex topological vector
        space.
    \item Every class $[f] \in \Aph$ has a unique real-analytic
        representative $f \in C^\omega(\mathbb{C}^n)$. Therefore, we
        identify $\Aph$ with the corresponding subspace of
        $C^\omega(\mathbb{C}^n)$ from now on.
    \item Every function $f \in \Aph$ has a unique extension to a
        function $\hat{f} \in \HolAntiHol (\mathbb{C}^n \times
        \mathbb{C}^n)$, i.e. holomorphic in the first and
        anti-holomorphic in the second argument, such that
        \begin{equation}
            \label{eq:fDeltahatf}
            f = \Delta^* \hat{f},
        \end{equation}
        where $\Delta: \mathbb{C}^n \ni z \mapsto (z,z) \in
        \mathbb{C}^n \times \mathbb{C}^n$ is the diagonal.
    \item Any $f \in C^\omega(\mathbb{C}^n)$ such that there exist
        constants $a, b, c > 0$ with
        \begin{equation}
            \label{eq:fabcEstimate}
            \left|
                \frac{\partial^{|R|+|S|} f}
                {\partial \cc{z}^R \partial z^S}(p)
            \right|
            \le c a^{|R|} b^{|S|}
        \end{equation}
        belongs to $\Aph$. In particular $\mathbb{C}[z, \cc{z}]
        \subseteq \Aph$.
    \end{enumerate}
\end{theorem}
\begin{proof}
    The first part is clear. For the second, we consider
    \[
    \normph[1,1,0,0]{f}
    \ge
    \frac{1}{R!} (\normph[0,0,R,0]{f})^2
    \ge
    \frac{(2\hbar)^{2|R|+2|S|}}{R!(S!)^2} 
    \left|
        \frac{\partial^{|R|+|S|} f}
        {\partial z^R \partial \cc{z}^S}(p) 
    \right|^4
    \]
    for all $R$, $S$. Thus we obtain
    \begin{equation}
        \label{eq:NormEinsTaylorCoeff}
        \left|
            \frac{\partial^{|R|+|S|} f}
            {\partial z^R \partial \cc{z}^S}(p) 
        \right|
        \le
        \normph[1,1,0,0]{f}
        \frac{\sqrt[4]{R!}\sqrt{S!}}{\sqrt{(2\hbar)^{|R| + |S|}}}.
    \end{equation}
    But this implies that the series
    \[
    \hat{f}(z, \cc{w})
    =
    \sum_{R, S=0}^\infty \frac{1}{R!S!}
    \frac{\partial^{|R|+|S|} f}
    {\partial z^R \partial \cc{z}^S}(p) 
    (z-p)^R(\cc{w} - \cc{p})^S
    \]
    converges for all $z, w \in \mathbb{C}^n$. Thus $\hat{f} \in
    \HolAntiHol(\mathbb{C}^n\times\mathbb{C}^n)$ and clearly
    $\Delta^*\hat{f}$ is in the same equivalence class as $f$. This
    shows the second and third part. Now assume $f \in
    C^\omega(\mathbb{C}^n)$ satisfies \eqref{eq:fabcEstimate}. Then
    \[
    \hph[0,0,R,S]{f}
    \le
    c^2 \E^{2\hbar b^2 n} (2\hbar a^2)^{|R|} (2 \hbar b^2)^{|S|},
    \]
    whence there are constants $c_{0,0} = c^2 \E^{2\hbar b^2 n}$,
    $a_{0,0} = 2\hbar a^2$ and $b_{0,0} = 2\hbar b$ such that
    \[
    \hph[0,0,R,S]{f}
    \le c_{0,0} \; a_{0,0}^{|R|} \; b_{0,0}^{|S|}.
    \]
    We claim that for all $m$, $\ell$ there are constants
    $a_{m,\ell}$, $b_{m, \ell}$ and $c_{m, \ell}$ such that
    \[
    \hph[m,\ell,R,S]{f}
    \le
    c_{m,\ell} \; a_{m,\ell}^{|R|} \; b_{m,\ell}^{|S|}.
    \]
    Indeed, a recursive argument shows that 
    \[
    c_{m,\ell} 
    = c_{m-1,\ell/2}^2 \E^{n(1 + b_{m-1, \ell/2})^2},
    \quad
    a_{m,\ell} = (1 + a_{m-1, \ell/2})^2,
    \quad
    b_{m,\ell} = (1 + b_{m-1, \ell/2})^2
    \]
    for even $\ell$ and
    \[
    c_{m,\ell} 
    = c_{m-1,(\ell-1)/2}^2 \E^{n(1 + a_{m-1, (\ell-1)/2})^2},
    \quad
    a_{m,\ell} = (1 + b_{m-1, (\ell-1)/2})^2,
    \quad
    b_{m,\ell} = (1 + a_{m-1, (\ell-1)/2})^2
    \]
    for odd $\ell$ will do the job. But then all seminorms of $f$ are
    finite.
\end{proof}

Since the polynomials are in $\Aph$ we shall make intense use of them.
The next proposition gives a first hint on the continuity of the Wick
product. Here and in the following we shall use the notation
$\starwickh$ for the Wick product with $\lambda$ being replaced by
$\hbar$.
\begin{proposition}
    \label{proposition:fWickPolynomialCont}
    Let $f \in C^\infty(\mathbb{C}^n)$ and let $g \in \mathbb{C}[z,
    \cc{z}]$ be a polynomial.
    \begin{enumerate}
    \item $f \starwickh g$ is a finite sum and thus a well-defined
        smooth function.
    \item $\normph{\cc{f} \starwickh g} \le
        \normph[m+1,2^m+\ell,R,S]{f} \normph[m+1,\ell,R,S]{g}$.
    \end{enumerate}
\end{proposition}
\begin{proof}
    The first part is clear from the explicit form of $\starwick$.
    For the second part we first have by the Cauchy-Schwarz inequality
    \[
    \begin{split}
        \left| \cc{f} \starwickh g \right|^2(p)
        &=
        \left| 
            \sum_{N=0}^{\textrm{finite}}
            \frac{(2\hbar)^{|N|}}{N!}
            \; \frac{\partial^{|N|} \cc{f}}{\partial z^N}(p)
            \; \frac{\partial^{|N|} g}{\partial \cc{z}^N}(p)
        \right|^2 \\
        &\le
        \left(
            \sum_{N=0}^\infty 
            \frac{(2\hbar)^{|N|}}{N!}
            \left|
                \frac{\partial^{|N|} \cc{f}}{\partial z^N}(p)
            \right|^2
        \right)
        \left(
            \sum_{N=0}^\infty 
            \frac{(2\hbar)^{|N|}}{N!}
            \left|
                \frac{\partial^{|N|} g}{\partial \cc{z}^N}(p)
            \right|^2
        \right) \\
        &=
        \hph[0,0,0,0]{f} \hph[0,0,0,0]{g}.
    \end{split}
    \tag{$*$}
    \]
    Now partial derivatives are still \emph{derivations} of
    $\starwick$ and hence of $\starwickh$ if one of the functions is a
    polynomial. This allows the following computation
    \[
    \begin{split}
        &\hph[0,0,R,S]{\cc{f} \starwickh g}\\
        &=
        \sum_{N=0}^\infty 
        \frac{(2\hbar)^{|N|+|R|+|S|}}{N!}
        \left|
            \frac{\partial^{|N|+|R|+|S|}}
            {\partial\cc{z}^{N+S}\partial z^R}
            \left(\cc{f} \starwickh g\right)
        \right|^2(p) \\
        &=
        \sum_{N=0}^\infty 
        \frac{(2\hbar)^{|N|+|R|+|S|}}{N!}
        \left|
            \sum_{I=0}^R\sum_{J=0}^{N+S}
            \Binom{R}{I}\Binom{N+S}{J}
            \frac{\partial^{|I|+|J|} \cc{f}}
            {\partial z^I \partial \cc{z}^J}
            \starwickh
            \frac{\partial^{|N+S-J|+|R-I|} g}
            {\partial z^{R-I} \partial \cc{z}^{N+S-J}}
        \right|^2(p) \\
        &\le
        \sum_{N=0}^\infty 
        \frac{(2\hbar)^{|N|+|R|+|S|}}{N!}
        \left(
            \sum_{I=0}^R\sum_{J=0}^{N+S}
            \Binom{R}{I}\Binom{N+S}{J}
            \left|
                \frac{\partial^{|I|+|J|} \cc{f}}
                {\partial z^I \partial \cc{z}^J}
                \starwickh
                \frac{\partial^{|N+S-J|+|R-I|} g}
                {\partial z^{R-I} \partial \cc{z}^{N+S-J}}
            \right|
        \right)^2(p) \\
        &\stackrel{(*)}{\le}
        \sum_{N=0}^\infty
        \frac{(2\hbar)^{|N|+|R|+|S|}}{N!}
        \left(
            \sum_{I=0}^R\sum_{J=0}^{N+S}
            \Binom{R}{I}\Binom{N+S}{J}
            \sqrt{
              \hph[0,0,0,0]{
                \frac{\partial^{|I|+|J|} f}
                {\partial \cc{z}^I \partial z^J}
              }
            }
            \sqrt{
              \hph[0,0,0,0]{
                \frac{\partial^{|N+S-J|+|R-I|} g}
                {\partial z^{R-I} \partial \cc{z}^{N+S-J}}
              }
            }
        \right)^2 \\
        &=
        \sum_{N=0}^\infty
        \frac{1}{N!}
        \left(
            \sum_{I=0}^R\sum_{J=0}^{N+S}
            \Binom{R}{I}\Binom{N+S}{J}
            \sqrt{\hph[0,0,J,I]{f}}
            \sqrt{\hph[0,0,R-I,N+S-J]{g}}
        \right)^2 \\
        &\le
        \sum_{N=0}^\infty
        \frac{1}{N!}
        \left(
            \sum_{I=0}^R\sum_{J=0}^{N+S}
            \Binom{R}{I}\Binom{N+S}{J}
            \hph[0,0,J,I]{f}
        \right)
        \left(
            \sum_{I=0}^R\sum_{J=0}^{N+S}
            \Binom{R}{I}\Binom{N+S}{J}
            \hph[0,0,I,J]{g}
        \right) \\
        &\le
        \sqrt{
          \sum_{N=0}^\infty
          \frac{1}{N!}
          \left(
              \sum_{I=0}^R\sum_{J=0}^{N+S}
              \Binom{R}{I}\Binom{N+S}{J}
              \hph[0,0,J,I]{f}
          \right)^2
        } 
        \sqrt{
          \sum_{N=0}^\infty
          \frac{1}{N!}
          \left(
              \sum_{I=0}^R\sum_{J=0}^{N+S}
              \Binom{R}{I}\Binom{N+S}{J}
              \hph[0,0,I,J]{g}
          \right)^2
        }\\
        &=
        \sqrt{\hph[1,1,R,S]{f}} \sqrt{\hph[1,0,R,S]{g}},
    \end{split}
    \]
    which proves the second part for the case $m=0$. Note the
    necessity that one function (we have chosen $g$) is polynomial
    since otherwise the `function' $\cc{f} \starwickh g$ is a priori
    not defined as a smooth function. The general case is now obtained
    by a straightforward induction on $m$ only using the
    Cauchy-Schwarz inequality.
\end{proof}
\begin{corollary}
    \label{corollary:AphStableUndercc}
    The pointwise complex conjugation is a continuous map $\Aph
    \longrightarrow \Aph$. We have
    \begin{equation}
        \label{eq:normccf}
        \normph{\cc{f}} \le 
        \normph[m+1,\ell,R,S]{1} \normph[m+1, 2^m + \ell, R, S]{f}.
    \end{equation}
\end{corollary}

Note however that the seminorms themselves are not invariant under
complex conjugation $f \leftrightarrow \cc{f}$, though the complex
conjugation is continuous.

We also note that the second part of the proposition already shows
some nice continuity properties of the Wick star product, at least if
one function is a polynomial. Note however, that the above argument
will not extend to arbitrary $f, g$ whence we shall need another
route.
\begin{theorem}
    \label{theorem:TaylorConverges}
    The polynomials $\mathbb{C}[z, \cc{z}] \subseteq \Aph$ are
    dense. More specifically, the Taylor expansion
    \begin{equation}
        \label{eq:frealanalytic}
        f(z, \cc{z})
        = \sum_{I, J = 0}^\infty \frac{1}{I!J!}
        \frac{\partial^{|I|+|J|}f}
        {\partial z^I \partial \cc{z}^J}(p)
        (z-p)^I (\cc{z} - \cc{p})^J
    \end{equation}
    of $f \in \Aph \subseteq C^\omega(\mathbb{C}^n)$ converges
    unconditionally to $f$ with respect to the topology of $\Aph$.  In
    particular, the truncated Taylor polynomials
    \begin{equation}
        \label{eq:fMNDef}
        f^{(N, M)}(z, \cc{z})
        = \sum_{I=0}^N \sum_{J = 0}^M 
        \frac{1}{I!J!}
        \frac{\partial^{|I|+|J|}f}
        {\partial z^I \partial \cc{z}^J}(p)
        (z-p)^I (\cc{z} - \cc{p})^J
    \end{equation}
    converge unconditionally to $f$ in the topology of $\Aph$.
\end{theorem}
\begin{proof}
    Clearly, we only have to show the later statement. First we
    rewrite the seminorms $\normph{\cdot}$ in the following
    `measure-theoretic' way
    \[
    \hph{f} = 
    \sum_{I_1, J_1 = 0}^\infty
    \cdots
    \sum_{I_s, J_s = 0}^\infty
    \mu^{m,\ell,R,S,\hbar}_{I_1, \cdots I_s, J_1, \cdots J_s}
    \left|
        \frac{\partial^{|I_1|+|J_1|}f}
        {\partial z^{I_1} \partial \cc{z}^{J_1}}(p)
    \right|^2
    \cdots
    \left|
        \frac{\partial^{|I_s|+|J_s|}f}
        {\partial z^{I_s} \partial \cc{z}^{J_s}}(p)
    \right|^2,
    \]
    where $\mu^{m,\ell,R,S,\hbar}_{I_1, \cdots I_s, J_1, \cdots, J_s}
    \ge 0$ are numerical constants not depending on $p$, and $s =
    2^m$.  This can be seen by induction easily. The concrete form of
    the coefficients $\mu^{m,\ell,R,S,\hbar}_{I_1, \cdots, I_s, J_1,
      \cdots, J_s}$ is not important for the following argument. Now
    we define for $f \in C^\infty(\mathbb{C}^n)$ a non-negative
    function $\phi_p(f): \mathbb{N}^{2sn} \longrightarrow [0, \infty)$
    by
    \[
    \phi_p(f)(I_1, \ldots, I_s, J_1, \ldots, J_s)
    =
    \left|
        \frac{\partial^{|I_1|+|J_1|}f}
        {\partial z^{I_1} \partial \cc{z}^{J_1}}(p)
    \right|^2
    \cdots
    \left|
        \frac{\partial^{|I_s|+|J_s|}f}
        {\partial z^{I_s} \partial \cc{z}^{J_s}}(p)
    \right|^2.
    \]
    Then we can interpret $\hph{f}$ as the `integral' of $\phi_p(f)$
    over $\mathbb{N}^{2sn}$ with respect to the weighted counting
    measure $\D\mu^{m,\ell,R,S,\hbar}$ determined by the coefficients
    $\mu^{m,\ell,R,S,\hbar}_{I_1, \cdots, I_s, J_1, \cdots, J_s}$,
    i.e.
    \[
    \hph{f} = \int_{\mathbb{N}^{2sn}} \phi_p(f) \D\mu^{m,\ell,R,S,\hbar}.
    \]
    Now let $\mathsf{K}, \mathsf{L} \subseteq \mathbb{N}^n$ be finite
    subsets and define the polynomial
    \[
    f^{(\mathsf{K}, \mathsf{L})} (z, \cc{z})
    =
    \sum_{I \in \mathsf{K}} \sum_{J \in \mathsf{J}}
    \frac{1}{I!J!} 
    \frac{\partial^{|I|+|J|f}}{\partial z^I \partial \cc{z}^J}(p) 
    (z-p)^I(\cc{z} - \cc{p})^J.
    \]
    Then we clearly have
    \begin{equation}
        \label{eq:ffKL}
        \phi_p\left(f- f^{(\mathsf{K}, \mathsf{L})} \right)
        (I_1, \ldots, I_s, J_1, \ldots, J_s)
        =
        \begin{cases}
            0
            & \textrm{if} \quad 
            I_1, \ldots, I_s \in \mathsf{K}, J_1, \ldots, J_s \in \mathsf{L}\\
            \phi_p(f)
            & \textrm{else},
        \end{cases}
    \end{equation}
    Thus when $\mathsf{K}, \mathsf{L}$ exhaust $\mathbb{N}^n$, the
    function $\phi_p(f - f^{(\mathsf{K}, \mathsf{L})})$ converges
    pointwise and monotonically to zero, i.e. for $\mathsf{K}
    \subseteq \mathsf{K}'$ and $\mathsf{L} \subseteq \mathsf{L}'$ we
    have
    \[
    \phi_p\left(f- f^{(\mathsf{K}, \mathsf{L})} \right)
    \ge
    \phi_p\left(f- f^{(\mathsf{K}', \mathsf{L}')} \right),
    \]
    and for all $I_1, \ldots, I_s, J_1, \ldots, J_s$
    \[
    \lim_{\mathsf{K}, \mathsf{L} \to \mathbb{N}^n}
    \phi_p\left(f- f^{(\mathsf{K}, \mathsf{L})} \right)
    (I_1, \ldots, I_s, J_1, \ldots, J_s)
    = 0
    \]
    in these sense of net convergence for the net of finite subsets of
    $\mathbb{N}^n$. Now an order of summation in
    \eqref{eq:frealanalytic} corresponds to a strictly increasing
    sequence $\mathsf{K}_i \times \mathsf{L}_i \subseteq \mathbb{N}^n
    \times \mathbb{N}^n$ which exhausts $\mathbb{N}^n \times
    \mathbb{N}^n$.  Then
    \[
    \lim_{i\to \infty}
    \hph{f - f^{(\mathsf{K}_i, \mathsf{L_i})}}
    =
    \lim_{i \to \infty}
    \int \phi_p\left(f - f^{(\mathsf{K}_i, \mathsf{L_i})}\right)
    \D\mu^{m,\ell,R,S,\hbar} = 0
    \]
    by dominated convergence. But this is equivalent to the
    unconditional convergence of \eqref{eq:frealanalytic}, see also
    \cite[Sect.~14.6, Thm.~1]{jarchow:1981a}.
\end{proof}

We now come to the main result of this section:
\begin{theorem}
    \label{theorem:AphFrechet}
    The locally convex topology of $\Aph$ is complete, i.e. $\Aph$ is
    a Fr\'echet space.
\end{theorem}
\begin{proof}
    Since the topology is determined by countably many seminorms we
    only have to consider Cauchy sequences and not Cauchy nets. Thus
    let $f_i \in \Aph$ be a Cauchy sequence, i.e. for all seminorms
    $\normph{\cdot}$ and all $\epsilon > 0$ we find a
    $K(m,\ell,R,S,\epsilon)$ such that for $i, j \ge K$ we have
    \[
    \normph{f_i - f_j} < \epsilon.
    \]
    We first evaluate this for $m = 1, \ell = 1, R, S = 0$. Let
    \[
    f_i(z, \cc{z}) = \sum_{I, J = 0}^\infty \frac{1}{I!J!}
    a^{(i)}_{IJ} (z - p)^I (\cc{z} - \cc{p})^J
    \]
    be the Taylor expansion of $f_i$ then from
    \eqref{eq:NormEinsTaylorCoeff} we see that the Taylor coefficients
    $a^{(i)}_{IJ}$ form a Cauchy sequence for each $I, J$. Denote
    their limit by
    \[
    a_{IJ} = \lim_{i \to \infty} a^{(i)}_{IJ}
    \tag{$*$}
    \]
    and define
    \[
    f(z, \cc{z}) = 
    \sum_{I, J = 0}^\infty \frac{1}{I!J!} a_{IJ} 
    (z - p)^I (\cc{z} - \cc{p})^J.
    \]
    Then we want to show $f \in \Aph$ and $f_i \to f$. To this end we
    first choose a smooth function $\tilde{f} \in
    C^\infty(\mathbb{C}^n)$ with Taylor coefficients at $p$ given by
    ($*$), which is possible thanks to the Borel Lemma. Since $f_i$ is
    a Cauchy sequence with respect to $\normph{\cdot}$ the sequence of
    seminorms $\normph{f_i}$ stays bounded as $i \to 0$. Thus we can
    again use the measure-theoretic point of view and write with the
    notation from the previous proof
    \[
    \begin{split}
        \hph{\tilde{f}} 
        &= 
        \int_{\mathbb{N}^{2sn}}
        \phi(\tilde{f}) \D \mu^{m,\ell,R,S,\hbar} \\
        &= 
        \int_{\mathbb{N}^{2sn}}
        \lim_{i\to \infty} \phi(f_i) \D \mu^{m,\ell,R,S,\hbar} \\
        &=
        \int_{\mathbb{N}^{2sn}}
        \liminf_i \phi(f_i) \D \mu^{m,\ell,R,S,\hbar} \\
        &\le
        \liminf_i \int_{\mathbb{N}^{2sn}}
        \phi(f_i) \D \mu^{m,\ell,R,S,\hbar} \\
        &\le
        \sup_i (\normph{f_i})^{2^{m+1}} < \infty,
    \end{split}
    \]
    by Fatou's Lemma. Thus $\tilde{f} \in \tAph$ and hence $f$ as in
    ($*$) is the unique real-analytic representative in $\Aph$. Next
    we compute using \eqref{eq:fMNDef}
    \[
    \begin{split}
        &\normph{f-f_i} \\
        &\le
        \normph{(f-f_i) - (f-f_i)^{(N,M)}}
        +
        \normph{f^{(N,M)}-f_j^{(N,M)}}
        +
        \normph{f_j^{(N,M)} - f_i^{(N,M)}} \\
        &\le
        \normph{(f-f_i) - (f-f_i)^{(N,M)}}
        +
        \normph{f^{(N,M)}-f_j^{(N,M)}}
        +
        \normph{f_j - f_i}.
    \end{split}
    \]
    Now we fix $\epsilon > 0$ and $K$ such that the last term is
    smaller than $\epsilon/3$ for $i,j > K$. For such an $i$ we fix
    $N, M$ such that the first term is smaller $\epsilon/3$ thanks to
    Theorem~\ref{theorem:TaylorConverges}. Finally, for this choice of
    $N, M$ we can find $j$ large enough that the second term is
    smaller $\epsilon/3$ since we have two polynomials of \emph{fixed}
    degree $(N, M)$ whose coefficients converge. This finally proves
    $f_i \to f$ with respect to $\normph{\cdot}$.
\end{proof}

%
%

\section{The continuity of $\starwick$}
\label{sec:ContinuityWick}

From Proposition~\ref{proposition:fWickPolynomialCont} we know that on
the subspace $\mathbb{C}[z, \cc{z}] \subseteq \Aph$ the Wick star
product $\starwickh$ is well-defined and continuous with respect to
the topology of $\Aph$. Since on the other hand $\mathbb{C}[z,
\cc{z}]$ is a dense subspace by Theorem~\ref{theorem:TaylorConverges}
the Wick star product extends uniquely to a continuous product on
$\Aph$ which thereby becomes a Fr\'echet algebra. However, from this
abstract extension we cannot yet conclude whether the \emph{formula}
\eqref{eq:WickMultiindex} with $\lambda$ being replaced by $\hbar$ is
still true. Thus we need an additional argument.

Let $\ell \in \{0, \ldots, 2^m -1\}$ be written as $\ell = \ell_{m-1}
2^{m-1} + \cdots + \ell_1 2 + \ell_0$ with $\ell_{m-1}, \ldots, \ell_0
\in \{0,1\}$. Then we define $\epsilon_\ell = (-1)^{\ell_{m-1} +
  \cdots + \ell_0}$. With this notation we can prove the following
continuity property of the partial derivatives:
\begin{proposition}
    \label{proposition:PartialDerivCont}
    Let $f \in C^\infty(\mathbb{C}^n)$ then we have
    \begin{equation}
        \label{eq:PartDerivCont}
        \sqrt{(2\hbar)^{|I|+|J|}}
        \normph{
          \frac{\partial^{|I|+|J|}f}
          {\partial z^I \partial \cc{z}^J}
        }
        \le
        \begin{cases}
            \normph[m,\ell,R+I,S+J]{f} 
            & \textrm{for} \quad \epsilon_\ell = +1 \\
            \normph[m,\ell,R+J,S+I]{f}
            & \textrm{for} \quad \epsilon_\ell = -1.
        \end{cases}
    \end{equation}
\end{proposition}
\begin{proof}
    Clearly, for $m=0$ (and hence $\ell = 0$) we even have the
    equality
    \[
    \sqrt{(2\hbar)^{|I|+|J|}}
    \normph[0,0,R,S]{
      \frac{\partial^{|I|+|J|}f}
      {\partial z^I \partial \cc{z}^J}
    }
    =
    \normph[0,0,R+I,S+J]{f}
    \]
    by the very definition of $\normph[0,0,R,S]{\cdot}$. Thus we prove
    the claim by induction on $m$. Let first be $\ell$ even and
    $\epsilon_\ell = +1$. Then $\epsilon_{\ell/2} = +1$ as well and we
    have by induction
    \[
    \begin{split}
        \hph{
          \sqrt{(2\hbar)^{|I|+|J|}}
          \frac{\partial^{|I|+|J|}f}
          {\partial z^I \partial \cc{z}^J}
        }
        &\le
        \sum_{N=0}^\infty \frac{1}{N!}
        \left(
            \sum_{K=0}^R\sum_{L=0}^{N+S}
            \Binom{R}{K}\Binom{N+S}{L}
            \hph[m-1,\ell/2,K+I, L+J]{f}
        \right)^2 \\
        &=
        \sum_{N=0}^\infty \frac{1}{N!}
        \left(
            \sum_{K=I}^{R+I}\sum_{L=J}^{N+S+J}
            \Binom{R}{K-I}\Binom{N+S}{L-J}
            \hph[m-1,\ell/2,K,L]{f}
        \right)^2 \\
        &\le
        \sum_{N=0}^\infty \frac{1}{N!}
        \left(
            \sum_{K=I}^{R+I}\sum_{L=J}^{N+S+J}
            \Binom{R+I}{K}\Binom{N+S+J}{L}
            \hph[m-1,\ell/2,K,L]{f}
        \right)^2 \\
        &\le
        \hph[m,\ell,R+I,S+J]{f}.
    \end{split}
    \]
    For $\epsilon_\ell = \epsilon_{\ell/2} = -1$ we get
    $\hph[m,\ell,R+J,S+I]{f}$ instead and the two cases with odd
    $\ell$ are analogous.
\end{proof}

From Theorem~\ref{theorem:TaylorConverges} we see that the polynomials
\begin{equation}
    \label{eq:zetapIJDef}
    \zeta_p^{IJ}(z, \cc{z}) = (z-p)^I (\cc{z} - \cc{p})^J
\end{equation}
form a countable unconditional topological basis for $\Aph$. The
proposition now implies that we even have a \emph{Schauder basis}:
\begin{corollary}
    \label{corollary:SchauderBasis}
    The polynomials $\{\zeta_p^{IJ}\}_{I, J \in \mathbb{N}^n}$ form an
    unconditional Schauder basis for $\Aph$.
\end{corollary}
\begin{proof}
    By Theorem~\ref{theorem:TaylorConverges} we have the unconditional
    convergence
    \[
    f = \sum_{I, J = 0}^\infty \frac{1}{I!J!} 
    \delta_p^{IJ}(f) \zeta_p^{IJ},
    \]
    and the $(I, J)$-th derivative $\delta_p^{IJ}$ of the
    $\delta$-functionals are continuous linear functionals on $\Aph$
    by Proposition~\ref{proposition:PartialDerivCont}. This implies
    the result, see e.g.~\cite[Sect.~14.2]{jarchow:1981a} for a
    definition of a Schauder basis.
\end{proof}

The next proposition shows that the pointwise product is continuous in
the topology of $\Aph$:
\begin{proposition}
    \label{proposition:PointwiseProductContinuous}
    Let $f, g \in C^\infty(\mathbb{C}^n)$. Then we have
    \begin{equation}
        \label{eq:NormphPointwiseProduct}
        \normph{fg}
        \le
        \normph[m+1,\ell,R,S]{f}\normph[m+1,\ell,R,S]{g},
    \end{equation}
    whence $\Aph$ is a Fr\'echet $^*$-algebra with respect to the
    pointwise product.
\end{proposition}
\begin{proof}
    First recall that from the explicit form of the Wick star product
    we obtain
    \[
    (2\hbar)^{|I|+|J|} 
    \left|
        \frac{\partial^{|I|+|J|}f}
        {\partial z^I \partial \cc{z}^J}
    \right|^2 (p)
    \le \hph[0,0,I,J]{f}.
    \tag{$*$}
    \]
    Now we first consider the case $m = 0$. Here we have by the
    Leibniz rule
    \[
    \begin{split}
        &\hph[0,0,R,S]{fg}\\
        &=
        \sum_{N=0}^\infty \frac{(2\hbar)^{|N|+|R|+|S|}}{N!}
        \left|
            \sum_{I=0}^R\sum_{J=0}^{N+S}
            \Binom{R}{I}\Binom{N+S}{J}
            \frac{\partial^{|I|+|J|}f}
            {\partial z^I \partial \cc{z}^J}
            \frac{\partial^{|R-I|+|N+S-J|}g}
            {\partial z^{R-I} \partial \cc{z}^{N+S-J}}
        \right|^2(p) \\
        &\le
        \sum_{N=0}^\infty \frac{(2\hbar)^{|N|+|R|+|S|}}{N!}
        \left(
            \sum_{I=0}^R\sum_{J=0}^{N+S}
            \Binom{R}{I}\Binom{N+S}{J}
            \left|
                \frac{\partial^{|I|+|J|}f}
                {\partial z^I \partial \cc{z}^J}(p)
            \right|
            \left|
                \frac{\partial^{|R-I|+|N+S-J|}g}
                {\partial z^{R-I} \partial \cc{z}^{N+S-J}}(p)
            \right|
        \right)^2 \\
        &\stackrel{(*)}{\le}
        \sum_{N=0}^\infty \frac{1}{N!}
        \left(
            \sum_{I=0}^R\sum_{J=0}^{N+S}
            \Binom{R}{I}\Binom{N+S}{J}
            \sqrt{\hph[0,0,I,J]{f}}
            \sqrt{\hph[0,0,R-I,N+S-J]{g}}
        \right)^2 \\
        &\le
        \sum_{N=0}^\infty \frac{1}{N!}
        \left(
            \sum_{I=0}^R\sum_{J=0}^{N+S}
            \Binom{R}{I}\Binom{N+S}{J}
            \hph[0,0,I,J]{f}
        \right)
        \left(
            \sum_{I=0}^R\sum_{J=0}^{N+S}
            \Binom{R}{I}\Binom{N+S}{J}
            \hph[0,0,I,J]{g}
        \right) \\
        &\le
        \sqrt{
          \sum_{N=0}^\infty \frac{1}{N!}
          \left(
              \sum_{I=0}^R\sum_{J=0}^{N+S}
              \Binom{R}{I}\Binom{N+S}{J}
              \hph[0,0,I,J]{f}
          \right)^2
        }\\
        &\qquad \times
        \sqrt{
          \sum_{N=0}^\infty \frac{1}{N!}
          \left(
              \sum_{I=0}^R\sum_{J=0}^{N+S}
              \Binom{R}{I}\Binom{N+S}{J}
              \hph[0,0,I,J]{g}
          \right)^2
        }\\
        &=
        \sqrt{\hph[1,0,R,S]{f}} \sqrt{\hph[1,0,R,S]{g}},
    \end{split}
    \]
    using twice the Cauchy Schwarz inequality in the last steps.  For
    $m \ge 1$ we proceed by a straightforward induction for the two
    cases of $\ell$ even and odd separately.
\end{proof}
\begin{corollary}
    \label{corollary:AphIsPoisson}
    $\Aph$ is a Fr\'echet-Poisson $^*$-algebra.
\end{corollary}
\begin{proof}
    Clearly, the canonical Poisson bracket is continuous as $\{f, g\}$
    is the sum of pointwise products of partial derivatives of $f$ and
    $g$.
\end{proof}

Combining the last two propositions we can finally show the continuity
of the Wick star product: we even show that the series
\eqref{eq:WickMultiindex} converges in the topology of $\Aph$ if we
replace $\lambda$ by any complex number $\alpha$:
\begin{theorem}
    \label{theorem:WickProductContinuous}
    The Wick star product
    \begin{equation}
        \label{eq:WickProducta}
        f \starwicka g
        =
        \sum_{N=0}^\infty \frac{(2\alpha)^{|N|}}{N!}
        \frac{\partial^{|N|}f}{\partial z^N}
        \frac{\partial^{|N|}g}{\partial \cc{z}^N}
    \end{equation}
    converges absolutely in the topology of $\Aph$ for all $\alpha \in
    \mathbb{C}$ and gives a continuous associative product. If $\alpha
    = \hbar > 0$ then $\Aph$ becomes a Fr\'echet $^*$-algebra with
    respect to $\starwickh$ as the product and the complex conjugation
    as the $^*$-involution.
\end{theorem}
\begin{proof}
    Let $f, g \in \Aph$ then we have for even $\ell$ and
    $\epsilon_\ell = +1$
    \[
    \begin{split}
        &\normph[m,\ell,0,0]{
          \sum_{N=0}^\infty \frac{(2\alpha)^{|N|}}{N!}
          \frac{\partial^{|N|}f}{\partial z^N}
          \frac{\partial^{|N|}g}{\partial \cc{z}^N}
        }\\
        &\le
        \sum_{N=0}^\infty
        \frac{|2\alpha|^{|N|}}{N!}
        \normph[m,\ell,0,0]{
          \frac{\partial^{|N|}f}{\partial z^N}
          \frac{\partial^{|N|}g}{\partial \cc{z}^N}
        }\\  
        &\stackrel{\textrm{Prop.~\ref{proposition:PointwiseProductContinuous}}}{\le}
        \sum_{N=0}^\infty
        \frac{|2\alpha|^{|N|}}{N!}
        \normph[m+1,\ell,0,0]{
          \frac{\partial^{|N|}f}{\partial z^N}
        }
        \normph[m+1,\ell,0,0]{
          \frac{\partial^{|N|}g}{\partial \cc{z}^N}
        }\\  
        &\stackrel{\textrm{Prop.~\ref{proposition:PartialDerivCont}}}{\le}
        \sum_{N=0}^\infty
        \frac{|\frac{\alpha}{\hbar}|^{|N|}}{N!}
        \normph[m+1,\ell,N,0]{f}
        \normph[m+1,\ell,0,N]{g}
        \\
        &\le
        \sqrt{
          \sum_{N=0}^\infty
          \frac{|\frac{\alpha}{\hbar}|^{|N|}}{N!}
          \left(\normph[m+1,\ell,N,0]{f}\right)^2
        }
        \sqrt{
          \sum_{N=0}^\infty
          \frac{|\frac{\alpha}{\hbar}|^{|N|}}{N!}
          \left(\normph[m+1,\ell,0,N]{g}\right)^2
        }\\
        &\stackrel{\textrm{Prop.~\ref{proposition:SeminormsAreSeminorms}}}{\le}
        \sqrt{
          \sum_{N=0}^\infty
          \frac{|\frac{\alpha}{\hbar}|^{|N|}}{N!}
          \left(
              \sqrt[2^{m+3}]{N!}
              \normph[m+2,2\ell+1,0,0]{f}
          \right)^2
        } 
        \sqrt{
          \sum_{N=0}^\infty
          \frac{|\frac{\alpha}{\hbar}|^{|N|}}{N!}
          \left(
              \sqrt[2^{m+3}]{N!}
              \normph[m+2,2\ell,0,0]{g}
          \right)^2
        }\\
        &= \underbrace{\left(
            \sum_{N=0}^\infty
            \frac{|\frac{\alpha}{\hbar}|^{|N|}}{N!}
            \sqrt[2^{m+2}]{N!}
        \right)}_{=c_m(\frac{\alpha}{\hbar})}
        \normph[m+2,2\ell+1,0,0]{f}
        \normph[m+2,2\ell,0,0]{g}.
    \end{split}
    \]
    Since $c_m(\frac{\alpha}{\hbar})$ converges for all $\alpha \in
    \mathbb{C}$ we have shown the convergence of
    \eqref{eq:WickProducta} with respect to
    $\normph[m,\ell,0,0]{\cdot}$ for even $\ell$ and $\epsilon_\ell =
    +1$. The other three cases are shown analogously. As the topology
    of $\Aph$ is already determined by the seminorms
    $\normph[m,\ell,0,0]{\cdot}$ the \emph{absolute} convergence in
    the topology of $\Aph$ follows. From the above estimate (and the
    analogous ones for odd $\ell$ etc.) one also obtains the
    continuity of $\starwicka$. If $\alpha = \hbar$ is real, then the
    complex conjugation is a $^*$-involution showing the last
    statement.
\end{proof}

\begin{corollary}
    \label{corollary:HolomorphicDeformation}
    The Wick star product $\starwicka$ is a holomorphic deformation of
    the pointwise product in the sense of
    \cite{pflaum.schottenloher:1998a}.
\end{corollary}
\begin{corollary}
    \label{corollary:NicerEstimateContinuityWick}
    For $f, g \in \Aph$ we have
    \begin{equation}
        \label{eq:NiceEstimateWick}
        \normph{\cc{f} \starwickh g}
        \le
        \normph[m+1,2^m+\ell,R,S]{f}
        \normph[m+1,\ell,R,S]{g}.
    \end{equation}
\end{corollary}
\begin{proof}
    This follows from
    Proposition~\ref{proposition:fWickPolynomialCont}, the density of
    $\mathbb{C}[z, \cc{z}]$ in $\Aph$ and the continuity of
    $\starwickh$.
\end{proof}

Though $\Aph$ becomes a Fr\'echet $^*$-algebra, the topology is
\emph{not} locally $m$-convex in the sense of \cite{michael:1952a},
see also \cite[App.~A]{pflaum.schottenloher:1998a}.  Recall that a
locally convex algebra is called locally $m$-convex if there exist a
set of seminorms $\norm{\cdot}_i$ defining the topology such that
$\norm{ab}_i \le \norm{a}_i \norm{b}_i$.  Such locally $m$-convex
algebras always have a holomorphic calculus which fails for $\Aph$:
\begin{example}
    \label{example:BadFunctions}
    We consider the entire function $f \in \mathcal{O}(\mathbb{C})$
    defined by
    \begin{equation}
        \label{eq:BadGuy}
        f(z) = \sum_{r=0}^\infty \frac{z^r}{\sqrt[4]{r!}}.
    \end{equation}
    Then it is easy to see that $f \starwickh \cc{f}$ evaluated at
    $z=0$ converges only for $\hbar =0$. Since clearly the
    $\delta$-functional $\delta_p : \Aph \longrightarrow \mathbb{C}$
    is continuous we conclude that $f \not\in \mathcal{A}_{0,\hbar}$.
    This shows that $\mathcal{A}_{0,\hbar}$ does \emph{not} allow a
    holomorphic functional calculus as $z \in \mathcal{A}_{0,\hbar}$
    and the $\starwickh$-Taylor expansion of $f$ would again coincide
    with $f$ since $\starwickh$-power of $z$ coincide with the
    corresponding pointwise powers. Analogous arguments apply also for
    $p \ne 0$ and higher dimensions $n \ge 1$.
\end{example}
\begin{corollary}
    \label{corollary:ANotLocallyMConvex}
    The topology of $\Aph$ is not locally $m$-convex with respect to
    the Wick star product $\starwicka$ for all $\alpha$ and there is
    no general holomorphic calculus for $\Aph$.
\end{corollary}

In the following we shall equip $\Aph$ always with the Wick star
product $\starwickh$.

\begin{remark}
    \label{remark:OMMYResult}
    At this point it would be interesting to compare our algebra to
    the construction obtained in
    \cite{omori.maeda.miyazaki.yoshioka:2000a}: Here the authors
    consider the Weyl-Moyal star product, which on the formal level is
    known to be equivalent to the Wick star product, and establish a
    convergence scheme to obtain a certain Fr\'echet algebra as the
    completion of the polynomials. However, their construction is
    rather different from ours whence it seems difficult to
    investigate whether the usual formal equivalence transformation
    survives the convergence conditions.
\end{remark}

%
%

\section{Translations and rescalings}
\label{sec:TranslationRescale}

We shall now discuss the dependence of $\Aph$ on the point $p \in
\mathbb{C}^n$ and on the value $\hbar > 0$. We start with the
dependence on the point $p$.

Let $\alpha \in \mathbb{C}^n$ then for $f \in \Aph$ we define
\begin{equation}
    \label{eq:HolTranslation}
    (\tau_\alpha f)(z, \cc{z}) 
    = \hat{f}(z+\alpha, \cc{z})
\end{equation}
and
\begin{equation}
    \label{eq:AntiHolTranslation}
    (\cc{\tau}_{\cc{\alpha}} f)(z, \cc{z}) 
    = \hat{f}(z, \cc{z}+\cc{\alpha}),
\end{equation}
which is well-defined according to
Theorem~\ref{theorem:RealAnalytic}, Part~\textit{iii.)}. Moreover, we
consider the functions
\begin{equation}
    \label{eq:EalphabetaFunctions}
    \mathsf{e}_{\cc{\alpha},\beta} (z, \cc{z}) 
    = \E^{\hbar\cc{\alpha}\beta} 
    \E^{\cc{\alpha}z + \beta \cc{z}},
\end{equation}
which are elements in $\Aph$ according to
Theorem~\ref{theorem:RealAnalytic}, Part~\textit{iv.)}. The following
lemma is a simple computation, the results of which are well-known in
the case of the formal Wick star product:
\begin{lemma}
    \label{lemma:ExpInduceTranslations}
    Let $\alpha, \beta, \gamma, \delta \in \mathbb{C}^n$. Then we have
    for all $f \in \Aph$:
    \begin{enumerate}
    \item $\mathsf{e}_{\cc{\alpha},\beta} \starwickh
        \mathsf{e}_{\cc{\gamma}, \delta} = \E^{\hbar(\cc{\alpha}\delta
          - \beta\cc{\gamma})} \mathsf{e}_{\cc{\alpha} + \cc{\gamma},
          \beta + \delta}$.
    \item $\mathsf{e}_{\cc{\alpha},\beta} \starwickh f =
        \mathsf{e}_{\cc{\alpha},\beta} \cc{\tau}_{2\hbar\cc{\alpha}}
        f$.
    \item $f \starwickh \mathsf{e}_{\cc{\alpha},\beta} =
        \mathsf{e}_{\cc{\alpha},\beta} \tau_{2\hbar\beta} f$.
    \item The maps $\tau_\alpha$ and $\cc{\tau}_{\cc{\alpha}}$ are
        continuous linear bijections
        \begin{equation}
            \label{eq:TauAph}
            \tau_\alpha, \cc{\tau}_{\cc{\alpha}}:
            \Aph \longrightarrow \Aph.
        \end{equation}
    \end{enumerate}
\end{lemma}
\begin{proof}
    The only non-trivial point here is that for $f \in \Aph$ we have
    \[
    \tau_\alpha f 
    = \sum_{N=0}^\infty \frac{\alpha^N}{N!}
    \frac{\partial^|N|f}{\partial z^N}
    \]
    and analogously for $\cc{\tau}_{\cc{\alpha}}$ since $f$ has a
    extension to $\hat{f} \in
    \HolAntiHol(\mathbb{C}^n\times\mathbb{C}^n)$.  Then the
    computations for the first three parts are folklore.  The last
    part follows from
    \[
    \tau_\alpha f 
    = \mathsf{e}_{0, -\frac{\alpha}{2\hbar}} 
    \left(
        f \starwickh \mathsf{e}_{0, \frac{\alpha}{2\hbar}}
    \right)
    \tag{$*$}
    \]
    and the continuity of the pointwise multiplication as well as of
    the continuity of $\starwickh$. The same argument applies for
    $\cc{\tau}_{\cc{\alpha}}$. From ($*$) one can easily work out
    explicit estimates for $\normph{\tau_\alpha f}$ and
    $\normph{\cc{\tau}_{\cc{\alpha}} f}$ using
    Proposition~\ref{proposition:PointwiseProductContinuous} and
    Corollary~\ref{corollary:NicerEstimateContinuityWick}.
\end{proof}
\begin{corollary}
    \label{corollary:ExponentialsUnitary}
    Let $\alpha, \beta \in \mathbb{C}^n$.
    \begin{enumerate}
    \item $\mathsf{e}_{\cc{\alpha},\beta} \in \Aph$ is invertible with
        respect to $\starwickh$ with inverse given by
        $\mathsf{e}_{-\cc{\alpha}, -\beta}$.
    \item $\mathsf{e}_{\cc{\alpha},\beta}$ is unitary iff $\alpha =
        -\beta$ since in general $\cc{\mathsf{e}_{\cc{\alpha},\beta}}
        = \mathsf{e}_{\cc{\beta},\alpha}$.
    \end{enumerate}
\end{corollary}

We introduce now the following notation. For $w \in \mathbb{C}^n$ we
denote the translation by $w$ by $T_w(z) = z+w$ whence we have the
corresponding pull-back on functions $(T_w^*f)(z) = f(z+w)$. Moreover,
we set
\begin{equation}
    \label{eq:Unitaryw}
    u_w = \mathsf{e}_{\frac{1}{2\hbar}\cc{w}, -\frac{1}{2\hbar} w}
    \in \Aph,
\end{equation}
which is a unitary element of $\Aph$ according to
Corollary~\ref{corollary:ExponentialsUnitary}.
\begin{proposition}
    \label{proposition:TranslationInnerAutos}
    The translation group $\mathbb{C}^n$ acts via pull-backs by
    continuous inner $^*$-auto\-mor\-phisms
    \begin{equation}
        \label{eq:TranslationInnerAuto}
        T^*_w = \Ad_{\mathrm{\scriptscriptstyle Wick}}(u_w)
    \end{equation}
    on $\Aph$. 
\end{proposition}
\begin{proof}
    This is a simple consequence of
    Lemma~\ref{lemma:ExpInduceTranslations} and the fact that $u_w$ is
    unitary.
\end{proof}
\begin{remark}
    \label{remark:NotSoEasy}
    We remark that on a heuristic level (or on polynomial functions
    only) the statement of this proposition is folklore. Note that it
    is clear, that for the \emph{formal} Wick star product the
    statement is \emph{wrong}: the translations are only \emph{outer}
    automorphisms as the elements $u_w$ are \emph{not} well-defined as
    formal series in $\lambda$. In fact, it was one of our main
    motivations to find a reasonably large algebra where the statement
    of the proposition is still true, extending the polynomials.
\end{remark}
The fact that the translations act by inner automorphisms will
immediately imply the following result:
\begin{theorem}
    \label{theorem:AphIsAhOnly}
    Let $p, p' \in \mathbb{C}^n$. Then
    \begin{equation}
        \label{eq:AphIsApprimeh}
        \Aph = \mathcal{A}_{p', \hbar}
    \end{equation}
    as Fr\'echet $^*$-algebras.
\end{theorem}
\begin{proof}
    Let $f \in \Aph$ be given and $w = p' - p$. Then we have on one
    hand
    \[
    \norm{f}^{p',\hbar}_{m,\ell,R,S} = 
    \normph{T^*_w f} < \infty
    \]
    since
    \[
    \frac{\partial^{|I|+|J|}f}
    {\partial z^I \partial \cc{z}^J}(p')
    =
    \frac{\partial^{|I|+|J|} (T^*_w f)}
    {\partial z^I \partial \cc{z}^J}(p),
    \]
    and in the definition of $\normph{\cdot}$ only the Taylor
    coefficients of $f$ at $p$ are used while the combinatorial
    coefficients in the construction are the same for all points $p'$,
    $p$. This shows $f \in \mathcal{A}_{p',\hbar}$ whence by symmetry
    the equality \eqref{eq:AphIsApprimeh} as vector spaces follows.
    On the other hand we have
    \[
    \norm{f}^{p',\hbar}_{m,\ell,R,S}
    =
    \normph{u_w \starwickh f \starwickh u_{-w}}
    \le
    \normph[m+1,2^m+\ell,R,S]{\cc{u_w}} 
    \normph[m+2,2^{m+1}+\ell,R,S]{f}
    \normph[m+2,\ell,R,S]{u_w}.
    \]
    This shows that the seminorms
    $\norm{\cdot}^{p',\hbar}_{m,\ell,R,S}$ can be estimated against
    the seminorms $\normph{\cdot}$ whence the topology of
    $\mathcal{A}_{p',\hbar}$ is coarser than the one of $\Aph$. By
    symmetry $p' \leftrightarrow p$ we see that they actually
    coincide. The algebraic structures are the same anyway whence the
    theorem is shown.
\end{proof}

Thus we see \emph{a posteriori} that the construction of the Fr\'echet
algebra $\Aph$ does not depend on our choice $p \in \mathbb{C}^n$.
Hence we can simply write $\Ah = \Aph$ from now on. Note however, that
the system of seminorms $\normph{\cdot}$ \emph{depends} on the choice
of $p$, only the induced topology is independent of $p$.

As in the formal case, \emph{all} $\delta$-functionals are positive:
\begin{corollary}
    \label{corollary:DeltaFunctionalsCont}
    All $\delta$-functionals
    \begin{equation}
        \label{eq:deltapCont}
        \delta_p: \Ah \longrightarrow \mathbb{C}
    \end{equation}
    are continuous positive linear functionals.
\end{corollary}
\begin{proof}
    The positivity is clear from \eqref{eq:WickProducta} with $\alpha
    = \hbar > 0$ and the continuity follows from the theorem as
    $|\delta_p(f)| \le \normph[0,0,0,0]{f}$.
\end{proof}
\begin{corollary}
    \label{corollary:AhTopFiner}
    The Fr\'echet topology of $\Ah$ is finer than the topology of
    pointwise convergence.
\end{corollary}

In the next step we want to analyze the continuity properties of the
group representation $w \mapsto T_w^*$ further. To this end we
consider the dependence of the elements $u_w$ on $w$. Since
\begin{equation}
    \label{eq:uwuvcocycle}
    u_w \starwickh u_v = \E^{-\frac{\I}{2\hbar} \IM(\cc{w}v)} u_{w+v}
\end{equation}
the map $w \mapsto u_w$ is \emph{not} a group morphism. Since
$\E^{-\frac{\I}{2\hbar}\IM(\cc{w}v)}$ is even a non-trivial group
cocycle we need to pass to the central extension of the translation
group $\mathbb{C}^n$ by this cocycle, i.e. to the Heisenberg group
$\mathsf{H}_n$. Here we use the convention that $\mathsf{H}_n =
\mathbb{C}^n \times \mathbb{R}$ with multiplication law
\begin{equation}
    \label{eq:HeisenbergProd}
    (w, c) \cdot (w', c') = (w + w', c + c' + \IM(\cc{w} w')).
\end{equation}
Then it follows that
\begin{equation}
    \label{eq:uwcDef}
    \mathsf{H}_n \ni (w, c) 
    \; \mapsto \;
    u_{(w,c)} = \E^{-\frac{\I}{2\hbar}c} u_w \in \mathrm{U}(\Ah)
\end{equation}
is a group morphism from $\mathsf{H}_n$ into the group of unitaries
$\mathrm{U}(\Ah)$ in $\Ah$. Since $\E^{-\frac{\I}{2\hbar} c}$ is
central, \eqref{eq:uwcDef} factors to the group morphism $w \mapsto
T^*_w$, i.e. we have
\begin{equation}
    \label{eq:AduwcAduw}
    \Ad_{\mathrm{\scriptscriptstyle Wick}}(u_{(w,c)}) 
    = \Ad_{\mathrm{\scriptscriptstyle Wick}}(u_w) = T^*_w
\end{equation}
for all $(w, c) \in \mathsf{H}_n$. The Lie algebra $\mathfrak{h}^n$ of
$\mathsf{H}_n$ can be identified with $\mathsf{H}_n$ via the
exponential map. Then the Lie bracket is given by $[(w,c), (w',c')] =
(0, 2\IM(\cc{w}w'))$. The next theorem shows that the group morphism
$\mathsf{H}_n \longrightarrow \mathrm{U}(\Ah)$ is analytic and induces
a Lie algebra morphism $\mathfrak{h}^n \longrightarrow \Ah$:
\begin{theorem}
    \label{theorem:HeisenbergAction}
    \begin{enumerate}
    \item The map $\mathfrak{h}^n \cong \mathsf{H}_n \ni (w,c) \mapsto
        u_{(w,c)}$ is analytic with respect to the topology of $\Ah$.
    \item The generator $J_{(w,c)}$ of the one-parameter group $t
        \mapsto u_{(tw, tc)}$ is given by
        \begin{equation}
            \label{eq:MomentumMapTranslations}
            J_{(w,c)}(z, \cc{z}) = 
            \left.\frac{\D}{\D t}\right|_{t=0}
            u_{(tw,tc)} = 
            - \frac{\I}{2\hbar} c 
            + \frac{1}{2\hbar}(\cc{w}z - w \cc{z}),
        \end{equation}
        and $\mathfrak{h}^n \ni (w,c) \mapsto J_{(w,c)} \in \Ah$ is a
        Lie algebra morphism where $\Ah$ is equipped with the
        $\starwickh$-commutator as Lie bracket.
    \item We have for all $k$
        \begin{equation}
            \label{eq:kthPowerJwc}
            \left.\frac{\D^k}{\D t^k}\right|_{t=0}
            u_{(tw,tc)} = 
            J_{(w,c)} 
            \underbrace{\starwickh \cdots \starwickh}_{k \;
              \textrm{times}}
            J_{(w,c)},
        \end{equation}
        and explicitly for $c = 0$
        \begin{equation}
            \label{eq:JwcPowerk}
            \left.\frac{\D^k}{\D t^k}\right|_{t=0}
            u_{(tw,0)} = 
            \sum_{\ell=0}^{\lfloor k/2 \rfloor} \frac{k!}{\ell!(k-2\ell)!}
            \left(-\frac{\cc{w}w}{4\hbar}\right)^\ell
            \left(\frac{\cc{w}z - w\cc{z}}{2\hbar}\right)^{k-2\ell}.
        \end{equation}
    \end{enumerate}
\end{theorem}
\begin{proof} 
    We have $u_{(w,c)}(z, \cc{z}) = \E^{-\frac{\I}{2\hbar}c -
      \frac{\cc{w}w}{4\hbar} + \frac{1}{2\hbar}(\cc{w}z - w \cc{z})}$.
    The first factors $\E^{-\frac{\I}{2\hbar}c}$ and
    $\E^{-\frac{\cc{w}w}{2\hbar}}$ are clearly analytic as they are
    analytic functions times a fixed element (the identity) in $\Ah$.
    For the remaining factor we see that the $(z, \cc{z})$-Taylor
    expansion, which converges unconditionally in the topology of
    $\Ah$ coincides up to numerical factors with the $(w,
    \cc{w})$-Taylor expansion of this function. Thus the $(w,
    \cc{w})$-Taylor expansion converges also in $\Ah$ unconditionally
    which shows the first part. The second part is a trivial
    computation. For the third part, the contribution of $c$ is not
    essential whence we discuss the case $c = 0$ only. A
    straightforward computation of the Taylor expansion in $t$ gives
    immediately \eqref{eq:JwcPowerk}. On the other hand, the $k$-th
    power of the linear function $J_w = J_{(w,0)}$ satisfies the
    recursion formula
    \[
    J_w^{\star k} 
    = J_w J_w^{\star(k-1)} 
    + (k-1) \frac{\cc{w}w}{2\hbar} J_w^{\star(k-2)},
    \]
    which follows easily from the fact that partial derivatives are
    derivations of $\starwickh$ and $\cc{w} \frac{\partial}{\partial
      \cc{z}} J_w = \frac{\cc{w}w}{2\hbar}$ is central. In a last step
    one shows that also the right hand side of \eqref{eq:JwcPowerk}
    satisfies this recursion with the same initial conditions for
    $k=0,1$.
\end{proof}
\begin{corollary}
    \label{corollary:StarExponential}
    Let $w \in \mathbb{C}^n$ then the $\starwickh$-exponential
    function
    \begin{equation}
        \label{eq:StarExpLinear}
        \Exp\left(t J_w\right)
        =
        \sum_{k=0}^\infty \frac{t^k}{k!} J_w \underbrace{\starwickh
          \cdots \starwickh}_{k \;\textrm{times}} J_w
        =
        u_{tw}
    \end{equation}
    converges unconditionally in the topology of $\Ah$ for all $t \in
    \mathbb{R}$.
\end{corollary}
\begin{remark}
    \label{remark:StarExp}
    Again, the importance of this corollary is not the explicit
    computation of the star exponential which is folklore. Instead, we
    have found a well-defined analytic framework where the formula
    actually converges inside an \emph{algebra} of functions. Note
    also, that in general we cannot expect such a convergence as
    $\Ah$ does not allow a holomorphic functional calculus in general.
\end{remark}

Let us now discuss the dependence on the parameter $\hbar$. From a
simple dimensional analysis we see that with our convention for the
Poisson bracket $\{z^k, \cc{z}^\ell\} = \frac{2}{\I} \delta^{k\ell}$
the coordinates have to have the physical dimension
$[\textrm{action}]^{\frac{1}{2}}$. Thus a rescaling of $\hbar$, which
physically is of course absurd, has to be reinterpreted as a rescaling
of the coordinates $(z, \cc{z})$: we are not changing the value
$\hbar$ but the unit system. On the other hand, from a purely
mathematical point of view we cannot distinguish these two
interpretations. As we do not have any additional absolute scale in
our approach, the corresponding algebras $\Ah$ and
$\mathcal{A}_{\hbar'}$ should be isomorphic in order to be physically
reasonable. The following theorem will show that this is indeed the
case.

We define for $\alpha > 0$ the diffeomorphism $R_\alpha: \mathbb{C}^n
\longrightarrow \mathbb{C}^n$
\begin{equation}
    \label{eq:Rescaler}
    R_\alpha (z) = \sqrt{\alpha} z,
\end{equation}
whose inverse is $R_{\frac{1}{\alpha}}$.
\begin{theorem}
    \label{theorem:Rescaling}
    The pull-back $R_\alpha^*$ induces an isomorphism of Fr\'echet
    $^*$-algebras
    \begin{equation}
        \label{eq:RescaleAalphahAh}
        R_\alpha^*: 
        \mathcal{A}_{\alpha\hbar} \longrightarrow \Ah.
    \end{equation}
    In particular, we have for all $f \in \mathcal{A}_{\alpha\hbar}$
    \begin{equation}
        \label{eq:Rescalefalpha}
        \norm{R_\alpha^* f}^{0, \hbar}_{m, \ell, R, S} = 
        \norm{f}^{0, \alpha\hbar}_{m, \ell, R, S}.
    \end{equation}
\end{theorem}
\begin{proof}
    Let $f \in \mathcal{A}_{\alpha\hbar}$ be given. Then we have
    \[
    \frac{\partial^{|I|+|J|} (R_\alpha^*f)}
    {\partial z^I \partial \cc{z}^J}
    =
    \sqrt{\alpha}^{|I|+|J|} R_\alpha^* 
    \left(
        \frac{\partial^{|I|+|J|} f}
        {\partial z^I \partial \cc{z}^J}
    \right)
    \]
    and thus
    \[
    \begin{split}
        \norm{R_\alpha^* f}^{0, \hbar}_{0, 0, R, S} 
        &= 
        \sum_{N=0}^\infty \frac{(2\hbar)^{|N|+|R|+|S|}}{N!}
        \left|
            \frac{\partial^{|N|+|R|+|S|} (R_\alpha^*f)}
            {\partial z^R \partial \cc{z}^{N+S}}(0)
        \right|^2 \\
        &=
        \sum_{N=0}^\infty \frac{(2\alpha\hbar)^{|N|+|R|+|S|}}{N!}
        \left|
            \frac{\partial^{|N|+|R|+|S|} f}
            {\partial z^R \partial \cc{z}^{N+S}}(0)
        \right|^2 \\
        &= \norm{f}^{0,\alpha\hbar}_{0,0,R,S}.
    \end{split}
    \]
    Since the higher seminorms are constructed in a purely
    combinatorial way out of $\norm{\cdot}^{0,\hbar}_{0,0,R,S}$, we
    can conclude \eqref{eq:Rescalefalpha}. This shows that
    $R_\alpha^*: \mathcal{A}_{\alpha\hbar} \longrightarrow \Ah$ is an
    isomorphism of Fr\'echet spaces as we can exchange the role of
    $\hbar$ and $\alpha\hbar$ by passing from $\alpha$ to
    $\frac{1}{\alpha}$. Then it is easy to see that $R_\alpha^*$ is an
    algebra morphism: we use the convergence of
    \eqref{eq:WickProducta} and the continuity of $R_\alpha^*$ to
    obtain
    \[
    \begin{split}
        R_\alpha^*(f
        \mathbin{\star^{\alpha\hbar}_{\mathrm{\scriptscriptstyle Wick}}}
        g)
        &= 
        R_\alpha^* \left(
            \sum_{N=0}^\infty \frac{(2\alpha\hbar)^{|N|}}{N!}
            \frac{\partial^{|N|}f}{\partial z^N}
            \frac{\partial^{|N|}g}{\partial \cc{z}^N}
        \right) \\
        &=
        \sum_{N=0}^\infty \frac{(2\hbar)^{|N|}}{N!}
        \frac{\partial^{|N|}(R_\alpha^*f)}{\partial z^N}
        \frac{\partial^{|N|}(R_\alpha^*g)}{\partial \cc{z}^N} \\
        &=
        R_\alpha^*f \starwickh R_\alpha^*g.
    \end{split}
    \]
    The compatibility with the complex conjugation $\cc{R_\alpha^* f}
    = R_\alpha^* \cc{f}$ is obvious.
\end{proof}
\begin{remark}
    \label{remark:CompareDifferentHbar}
    Of course we can also directly compare the seminorms for different
    values of $\hbar$. Clearly, one has
    \begin{equation}
        \label{eq:NullNormHbarHbarPrime}
        \norm{f}^{0,\hbar}_{0,0,R,S} \le \norm{f}^{0,\hbar'}_{0,0,R,S}
    \end{equation}
    for $\hbar \le \hbar'$ whence by induction
    \begin{equation}
        \label{eq:mellNormHbarHbarPrime}
        \norm{f}^{0,\hbar}_{m,\ell,R,S} \le \norm{f}^{0,\hbar'}_{m,\ell,R,S}
    \end{equation}
    as well. For $\hbar \le \hbar'$ this gives the inclusion
    \begin{equation}
        \label{AhbarAhbarprime}
        \mathcal{A}_{\hbar'} \subseteq \Ah.
    \end{equation}
\end{remark}

%
%

\section{The GNS construction and coherent states}
\label{sec:GNSCoherentStates}

We shall now discuss the GNS construction corresponding to the
positive $\delta$-functionals
\begin{equation}
    \label{eq:deltapAh}
    \delta_p: \Ah \longrightarrow \mathbb{C},
\end{equation}
now in the convergent situation.
\begin{proposition}
    \label{proposition:GelfandIdeal}
    Let $p \in \mathbb{C}^n$. Then the Gel'fand ideal of $\delta_p$ is
    given by
    \begin{equation}
        \label{eq:GelfandIdealdeltap}
        \mathcal{J}_p = 
        \left\{
            f \in \Ah
            \; \left| \;
                \forall I: \quad
                \frac{\partial^{|I|}f}{\partial \cc{z}^I}(p) = 0
            \right.
        \right\},
    \end{equation}
    and the GNS pre-Hilbert space $\mathfrak{D}_p = \Ah \big/
    \mathcal{J}_p$ is a Fr\'echet space in the natural way, where the
    topology of $\mathfrak{D}_p$ is determined by the seminorms
    \begin{equation}
        \label{eq:DpSeminorms}
        \normph{[f]} =
        \inf\left\{
            \normph{f + g} 
            \; \Big| \;
            g \in \mathcal{J}_p
        \right\}.
    \end{equation}
\end{proposition}
\begin{proof}
    The statement \eqref{eq:GelfandIdealdeltap} is obvious. Since
    $\delta_p$ is continuous, the Gel'fand ideal is a closed subspace
    of $\Ah$. Thus the quotient $\mathfrak{D}_p$ is again a Fr\'echet
    space by general arguments, see e.g.~\cite[Sect.~4.4,
    Prop.~1]{jarchow:1981a}.
\end{proof}

Since the translation group acts by inner $^*$-automorphisms we can
safely specialize to the case $p=0$ in a first step. Then we can
describe the quotient $\mathfrak{D}_0$ more explicitly:
\begin{theorem}
    \label{theorem:GNSConstruction}
    Let $f, g \in \Ah$.
    \begin{enumerate}
    \item The $\cc{z}$-Taylor expansion $\Psi: \Ah \longrightarrow
        \Ah$ defined by
        \begin{equation}
            \label{eq:PsiDef}
            \Psi: f
            \; \mapsto \;
            \left(
                z \mapsto 
                \Psi_f(z)
                =
                \sum_{I=0}^\infty
                \frac{1}{I!}
                \frac{\partial^{|I|}f}{\partial \cc{z}^I}(0) \cc{z}^I
            \right)
        \end{equation}
        is a continuous projection with
        \begin{equation}
            \label{eq:KerPsiGelfandIdeal}
            \ker \Psi = \mathcal{J}_0.
        \end{equation}
        In fact,
        \begin{equation}
            \label{eq:SeminormPsif}
            \normNullh{\Psi_f} \le
            \normNullh{f},
        \end{equation}
        where equality holds if and only if $f$ is anti-holomorphic.
    \item The quotient $\mathfrak{D}_0$ is canonically isomorphic as a
        Fr\'echet space to the image $\mathfrak{D} = \image\Psi$ of
        $\Psi$ via
        \begin{equation}
            \label{eq:DNullIsomorphD}
            \mathfrak{D}_0 \ni [f]
            \; \mapsto \;
            \Psi_f \in \mathfrak{D}.
        \end{equation}
        In particular,
        \begin{equation}
            \label{eq:DNullNormDNorm}
            \normNullh{[f]} = \normNullh{\Psi_f}.
        \end{equation}
    \item The space $\mathfrak{D}$ is a dense subspace of the
        Bargmann-Fock Hilbert space $\HBF$. The map
        \eqref{eq:DNullIsomorphD} is an isometry of pre-Hilbert spaces
        and the Fr\'echet topology of $\mathfrak{D}$ is finer than the
        topology induced from the Hilbert space $\HBF$.
    \item The GNS representation on $\mathfrak{D}_0$ induces via
        \eqref{eq:DNullIsomorphD} the Bargmann-Fock representation of
        $\Ah$ on $\mathfrak{D}$, explicitly given by
        \begin{equation}
            \label{eq:GNSRepExplicit}
            \pi(f) \Psi_g
            =
            \sum_{I=0}^\infty \frac{(2\hbar)^{|I|}}{I!}
            \left(
                \sum_{J=0}^\infty \frac{1}{J!}
                \frac{\partial^{|I|+|J|} f}
                {\partial z^I\partial \cc{z}^J}(0)
                \cc{z}^J
            \right)
            \frac{\partial^{|I|} \Psi_g}{\partial \cc{z}^I},
        \end{equation}
        where both series converge in the topology of $\mathfrak{D}$.
    \item The bilinear map $\Ah \times \mathfrak{D} \ni (f, \Psi_g)
        \mapsto \pi(f)\Psi_g \in \mathfrak{D}$ is continuous with
        respect to the Fr\'echet topologies of $\Ah$ and
        $\mathfrak{D}$, respectively.
    \end{enumerate}
\end{theorem}
\begin{proof}
    For the first part we have $\Psi_f(\cc{z}) = \hat{f}(0,\cc{z})$
    whence $\Psi_f$ is indeed a well-defined anti-holomorphic
    function.  Moreover, \eqref{eq:SeminormPsif} is clear from our
    consideration in the proof of
    Theorem~\ref{theorem:TaylorConverges} since in the seminorm of
    $\Psi_f$ simply less Taylor coefficients contribute compared to the
    corresponding seminorm of $f$. Thus $\Psi_f \in \Ah$ and $\Psi$ is
    continuous. From the explicit form of $\Psi$ the equality
    \eqref{eq:KerPsiGelfandIdeal} and $\Psi^2 = \Psi$ are obvious.
    For the second part we first notice that \eqref{eq:DNullIsomorphD}
    is well-defined and bijective since $f - \Psi_f \in \ker\Psi =
    \mathcal{J}_0$ using the fact that $\Psi$ is a projection.
    Moreover, \eqref{eq:DNullNormDNorm} follows directly from
    \eqref{eq:SeminormPsif} since obviously $\normNullh{[f]} \le
    \normNullh{\Psi_f}$. Since the inverse of
    \eqref{eq:DNullIsomorphD} is simply given by the well-defined and
    continuous map $\Psi_f \mapsto [f]$, we see that
    \eqref{eq:DNullIsomorphD} is indeed an isomorphism of Fr\'echet
    spaces. For the third part we compute explicitly the GNS inner
    product on $\mathfrak{D}_0$
    \[
    \begin{split}
        \SP{[f],[g]}_{\mathfrak{D}_0}
        &=
        \delta_0 (\cc{f} \starwickh g)\\
        &=
        \sum_{N=0}^\infty \frac{(2\hbar)^{|N|}}{N!}
        \; \frac{\partial^{|N|} \cc{f}}{\partial z^N}(0)
        \; \frac{\partial^{|N|} g}{\partial \cc{z}^N}(0)\\
        &=
        \sum_{N=0}^\infty \frac{(2\hbar)^{|N|}}{N!}
        \; \cc{\frac{\partial^{|N|} \Psi_f}{\partial \cc{z}^N}(0)}
        \; \frac{\partial^{|N|} \Psi_g}{\partial \cc{z}^N}(0)\\
        &=
        \SP{\Psi_f, \Psi_g}_{\mathrm{\scriptscriptstyle BF}},
    \end{split}
    \]
    whence \eqref{eq:DNullIsomorphD} is isometric. In particular
    $\mathfrak{D} \subseteq \HBF$ follows, as 
    \[
    \SP{\Psi_f, \Psi_f}_{\mathrm{\scriptscriptstyle BF}} 
    = (\cc{f} \starwickh f)(0) 
    = (\normNullh[0,0,0,0]{f})^2
    < \infty.
    \tag{$*$}
    \]
    Since $\mathbb{C}[\cc{z}] \subseteq \mathfrak{D} \subseteq \HBF$,
    the subspace $\mathfrak{D}$ is dense in $\HBF$. Moreover, since
    $\normNullh[0,0,0,0]{f} = \normNullh[0,0,0,0]{\Psi_f}$ the
    estimate ($*$) implies that the Fr\'echet topology of
    $\mathfrak{D}$ is finer than the topology induced from $\HBF$.
    This shows the third part. For the fourth part recall that the GNS
    representation is defined by $\varrho(f)[g] = [f \starwickh g]$
    which translates via \eqref{eq:DNullIsomorphD} into
    \[
    \begin{split}
        \pi(f)\Psi_g 
        &= \Psi_{f \starwickh g} \\
        &= \Psi_{
          \sum_{N=0}^\infty \frac{(2\hbar)^{|N|}}{N!}
          \frac{\partial^{|N|}f}{\partial z^N}
          \frac{\partial^{|N|}g}{\partial \cc{z}^N}
        }\\
        &\stackrel{(a)}{=}
        \sum_{N=0}^\infty \frac{(2\hbar)^{|N|}}{N!}
        \Psi_{
          \frac{\partial^{|N|}f}{\partial z^N}
          \frac{\partial^{|N|}g}{\partial \cc{z}^N}
        }\\
        &\stackrel{(b)}{=}
        \sum_{N=0}^\infty \frac{(2\hbar)^{|N|}}{N!}
        \Psi_{\frac{\partial^{|N|}f}{\partial z^N}}
        \Psi_{\frac{\partial^{|N|}g}{\partial \cc{z}^N}}
        \\
        &=
        \sum_{N=0}^\infty \frac{(2\hbar)^{|N|}}{N!}
        \left(
            \sum_{M=0}^\infty \frac{1}{M!}
            \frac{\partial^{|N|+|M|}f}
            {\partial z^N \partial \cc{z}^M}(0)
            \; \cc{z}^M
        \right)
        \frac{\partial^{|N|} \Psi_g}{\partial \cc{z}^N},
    \end{split}
    \]
    where $(a)$ holds since $\Psi$ is continuous and $(b)$ holds since
    obviously $\Psi$ is a homomorphism of the pointwise product.
    Finally, in the last step we have used that $\Psi$ commutes with
    derivatives in $\cc{z}$-direction. This shows the fourth part and
    the last part follows immediately from $\pi(f)\Psi_g = \Psi_{f
      \starwickh g} = \Psi_{f \starwickh \Psi_g}$ and the continuity
    of $\Psi$ and $\starwickh$.
\end{proof}
\begin{corollary}
    \label{corollary:BFRepInjective}
    The Bargmann-Fock representation of $\Ah$ is injective.
\end{corollary}
The following corollary is remarkable in so far as closed subspaces of
Fr\'echet spaces usually do not have complementary closed subspaces:
\begin{corollary}
    \label{corollary:AhDOplusJ}
    The algebra $\Ah$ decomposes into two complementary closed
    subspaces $\Ah = \mathfrak{D} \oplus \mathcal{J}_0$.
\end{corollary}
\begin{remark}
    \label{remark:AhIsCompletionofCCR}
    Since the Bargmann-Fock representation is injective and since
    \begin{equation}
        \label{eq:Annihilator}
        \pi(z^i) = 2\hbar \frac{\partial}{\partial \cc{z}^i} = a_i
    \end{equation}
    \begin{equation}
        \label{eq:Creator}
        \pi(\cc{z}^i) = \cc{z}^i = a_i^\dag
    \end{equation}
    are the annihilation and creation operators we find another
    interpretation of the algebra $\Ah$: it is a (rather large)
    completion of the polynomials in the creation and annihilation
    operators in a certain Fr\'echet topology. In particular, this
    completion contains the usual unitary generators of the Weyl
    algebra, i.e. the exponential functions of $a_i$ and $a_i^\dag$.
    The are given by $\pi(\mathsf{e}_{\cc{\alpha}, \beta})$ for suitable
    $\alpha, \beta \in \mathbb{C}^n$.
\end{remark}
\begin{remark}
    \label{remark:UnboundedTechnics}
    Since the Bargmann-Fock representation is a $^*$-representation of
    the $^*$-algebra $\Ah$ by (in general) unbounded operators with
    common domain $\mathfrak{D}$, one can investigate the resulting
    $O^*$-algebra of unbounded operators by techniques as developed in
    e.g. \cite{schmuedgen:1990a}. In particular, it would be
    interesting to find more concrete characterizations of the
    Fr\'echet topologies of $\mathfrak{D}$ and $\Ah$.
\end{remark}

We now discuss the action of the translation group $\mathbb{C}^n$ and
its central extension $\mathsf{H}_n$. The following statement is
obvious, the representation itself being well-known:
\begin{lemma}
    \label{lemma:HeisenbergRepUnitary}
    The map
    \begin{equation}
        \label{eq:HeisenbergRep}
        \mathsf{H}_n \ni (w, c) 
        \; \mapsto \;
        U_{(w, c)} = \pi(u_{(w,c)}) \in \mathrm{U}(\HBF)
    \end{equation}
    is a strongly continuous unitary representation of the Heisenberg
    group. Explicitly,
    \begin{equation}
        \label{eq:UwcExplicitly}
        (U_{(w, c)} \psi)(\cc{z})
        =
        \E^{-\frac{\I}{2\hbar} c - \frac{\cc{w}w}{2\hbar}}
        \E^{-\frac{w\cc{z}}{2\hbar}} \psi(\cc{z} + \cc{w})
    \end{equation}
    for $\psi \in \HBF$. It factors to a projective representation
    $U_w = U_{(w, 0)}$ of the translation group $\mathbb{C}^n$.
\end{lemma}
\begin{proof}
    Since $\pi$ is a $^*$-representation it follows that
    $\pi(u_{(w,c)})$ is a unitary operator defined on the dense domain
    $\mathfrak{D}$. Thus it extends to a unitary operator on $\HBF$.
    The group representation property is obvious from
    \eqref{eq:uwcDef}. The explicit formula is a simple consequence of
    \eqref{eq:GNSRepExplicit}. The fact that \eqref{eq:HeisenbergRep}
    is strongly continuous is well-known but can also be shown within
    our approach directly: let $\phi, \psi \in \mathfrak{D}$ then
    $g(w,c) = \SP{\psi, U_{(w,c)}\phi}_{\mathrm{\scriptscriptstyle
        BF}}$ is real-analytic since
    $\SP{\cdot,\cdot}_{\mathrm{\scriptscriptstyle BF}}$ and $\pi$ are
    continuous with respect to the Fr\'echet topology and $(w,c)
    \mapsto u_{(w,c)}$ is real-analytic according to
    Theorem~\ref{theorem:HeisenbergAction}. Since $g(0,0) = \SP{\psi,
      \phi}_{\mathrm{\scriptscriptstyle BF}}$ we see that on the dense
    domain $\mathfrak{D}$ the representation \eqref{eq:HeisenbergRep}
    is weakly continuous at the identity. But this implies that it is
    strongly continuous on the whole group $\mathsf{H}_n$ and on the
    whole Hilbert space $\HBF$. Since the contribution of $c$ is only
    an overall phase, the representation clearly factors to a
    projective representation of $\mathbb{C}^n$.
\end{proof}

Since we have a group action of $\mathsf{H}_n$ we can formulate now
the following covariance property of the Bargmann-Fock representation
which is an obvious consequence of the fact that the
$^*$-automorphisms are \emph{inner}.
\begin{theorem}
    \label{theorem:Covariance}
    The Bargmann-Fock representation is $\mathsf{H}_n$-covariant with
    respect to the action by $^*$-automorphisms on $\Ah$ and the
    action by unitaries on $\HBF$, i.e. we have
    \begin{equation}
        \label{eq:BFisCovariant}
        \pi(\Ad_{\mathrm{\scriptscriptstyle Wick}}(u_{(w,c)}) f)
        =
        U_{(w,c)} \pi(f) U_{(w,c)}^*
        =
        U_w \pi(f) U_w^*
    \end{equation}
    for all $(w, c) \in \mathsf{H}_n$ and $f \in \Ah$.
\end{theorem}

Since the GNS representation is cyclic with cyclic vector $\Psi_1 = 1$
we obtain \emph{coherent states} with respect to the representation of
$\mathsf{H}_n$. We define the \emph{coherent state vector}
$\psi_{(w,c)} \in \HBF$ by
\begin{equation}
    \label{eq:CoherentStateVector}
    \psi_{(w,c)} = U_{(w,c)}^{-1} \psi_1
\end{equation}
explicitly given by
\begin{equation}
    \label{eq:CoherentStateVectorExplicit}
    \psi_{(w,c)} (\cc{z})
    =
    \E^{\frac{\I}{2\hbar}c  + \frac{\cc{w}w}{4\hbar}}
    \E^{\frac{w\cc{z}}{2\hbar}}.
\end{equation}
From the covariance property \eqref{eq:BFisCovariant} we immediately
have the following characterization of the $\delta$-functionals at
arbitrary points in $\mathbb{C}^n$:
\begin{corollary}
    \label{corollary:CoherentStatesDeltap}
    The coherent state vectors give the $\delta_p$-functionals as
    expectation value functionals, i.e. we have
    \begin{equation}
        \label{eq:DeltapCoherent}
        \delta_w(f) = 
        \SP{\psi_{(w,c)}, 
          \pi(f) \psi_{(w,c)}}_{\mathrm{\scriptscriptstyle BF}}
    \end{equation}
    for all $f \in \Ah$ and $(w, c) \in \mathsf{H}_n$. In particular,
    the group action of $\mathsf{H}_n$ on the coherent state vectors
    factors through to a group action of the translation group
    $\mathbb{C}^n$ on the coherent states $\delta_w: \Ah
    \longrightarrow \mathbb{C}$.
\end{corollary}
\begin{remark}
    \label{remark:CoherentStates}
    This corollary gives finally the justification to view the
    $\delta$-functionals of the Wick star product algebra as coherent
    states with respect to the translation group. Though the explicit
    formula \eqref{eq:CoherentStateVectorExplicit} is folklore (it is
    just the Bergmann kernel) we would like to emphasize that in our
    approach the coherent states emerge out of properties of the
    observable algebra instead of more conventional approaches based
    on group actions on the state vectors in some Hilbert space, see
    e.g.  \cite{perelomov:1986a, yaffe:1982a} for more references. In
    this sense our approach supports the idea that the observable
    algebra is the more fundamental object in both, quantum and
    classical mechanics.
\end{remark}
\begin{remark}
    \label{remark:NotInTheFormalWay}
    We also note that the statement of
    Theorem~\ref{theorem:Covariance} as well as the
    Corollary~\ref{corollary:CoherentStatesDeltap} are \emph{not}
    possible for the formal Bargmann-Fock representation. Again, this
    was one of our main motivations to consider a suitable convergence
    scheme for the Wick star product.
\end{remark}

The result of Theorem~\ref{theorem:Covariance} and
Corollary~\ref{corollary:CoherentStatesDeltap} suggest the following
general definition of coherent states with respect to some symmetry
based on the \emph{observable algebra}:
\begin{definition}
    \label{definition:CoherentStates}
    Let $\mathcal{A}$ be a $^*$-algebra with unit $\Unit$ and let $G$
    be a group acting on $\mathcal{A}$ by $^*$-auto\-mor\-phisms
    $\Phi_g: \mathcal{A} \longrightarrow \mathcal{A}$. Let $\omega$ be
    a state of $\mathcal{A}$ such that the GNS representation is
    $G$-covariant, i.e. there exists a unitary (or more general:
    projectively unitary) representation $U$ of $G$ on the GNS
    pre-Hilbert space $\mathcal{H}_\omega$. Then the states $\omega_g$
    with
    \begin{equation}
        \label{eq:CoherentStatesDef}
        \omega_g(a) 
        = (\omega \circ \Phi_g)(a)
        = \SP{\psi_g, \pi(a) \psi_g},
    \end{equation}
    where $\psi_g = U_g^*\psi_\Unit \in \mathcal{H}_\omega$ are called
    coherent with respect to $G$.
\end{definition}

Clearly, in case of a projective representation only the coherent
states $\omega_g$ are well-defined, while for a unitary representation
also the coherent state vectors $\psi_g$ are well-defined.

With the Heisenberg group acting on $\Ah$ and the $\delta$-functional
we are in this situation: if the action is realized by \emph{inner}
$^*$-automorphisms then the representation is always covariant for a
(in general only projective) representation on the GNS pre-Hilbert
space. Note also, that for an \emph{invariant} state $\omega$ the GNS
representation is trivially covariant. In this case $\psi_g =
\psi_\Unit$ coincides with the vacuum vector for all $g \in G$. Thus
the interesting coherent states arise from non-invariant vacua such
that the GNS representation is nevertheless covariant. For a survey on
covariant $^*$-representation theory for $^*$-algebras over ordered
rings we refer to \cite{jansen.waldmann:2004a:pre}.

Let us now come to a further property of the subspace $\mathfrak{D}
\subseteq \HBF$. We have already seen that the action of
$\mathsf{H}_n$ leaves $\mathfrak{D}$ invariant.
\begin{theorem}
    \label{theorem:AnalyticVectors}
    The vectors in $\mathfrak{D}$ are analytic with respect to the
    unitary representation $U$ of $\mathsf{H}_n$.
\end{theorem}
\begin{proof}
    Let $\psi \in \mathfrak{D}$ be fixed then we have to show that
    $\mathsf{H}_n \ni (w, c) \mapsto U_{(w,c)}\psi$ is analytic with
    respect to the topology of $\HBF$. But this is simple since
    $U_{(w,c)} \psi = \pi(u_{(w,c)})\psi$ and $(w,c) \mapsto
    u_{(w,c)}$ is analytic in the topology of $\Ah$. Moreover, the
    \emph{bilinear} map $(f, \psi) \mapsto \pi(f)\psi$ is continuous
    in the topologies of $\Ah$ and $\mathfrak{D}$. Thus the map $(w,c)
    \mapsto U_{(w,c)}\psi$ is analytic with respect to the topology of
    $\mathfrak{D}$. Since by Theorem~\ref{theorem:GNSConstruction}
    this topology is finer than the one of $\HBF$, the proof is
    complete.
\end{proof}

Thus it would be interesting to know whether $\mathfrak{D}$ coincides
with the space of all analytic vectors. A positive answer would help
to understand the (still rather complicated) Fr\'echet topology of
$\mathfrak{D}$ and hence the one of $\Ah$.

Using the convergence of the star exponential \eqref{eq:StarExpLinear}
we even can specify the analyticity of the vectors in $\mathfrak{D}$
further. Since $\pi$ is continuous with respect to the topologies of
$\Ah$ and $\mathfrak{D}$, we have
\begin{equation}
    \label{eq:UwcStarExp}
    U_{(w,c)}\psi
    = \pi(u_{(w,c)})\psi
    = \pi\left(\Exp(J_{(w,c)})\right)\psi
    = \sum_{r=0}^\infty \frac{1}{r!} \pi\left(J_{(w,c)}\right)^r \psi
\end{equation}
with respect to the topology of $\mathfrak{D}$.  Again, since the
topology of $\mathfrak{D}$ is finer than the one of $\HBF$, the series
converges unconditionally also in the Hilbert space sense. As usual,
the contribution of $c$ is not essential.
\begin{corollary}
    \label{corollary:DanalyticCreatorAnnihilator}
    Let $\psi \in \mathfrak{D}$ and $w \in \mathbb{C}^n$. Then the
    series
    \begin{equation}
        \label{eq:UwConverges}
        U_w \psi 
        = \sum_{r=0}^\infty \frac{1}{r!}
        \left(\frac{1}{2\hbar}\right)^r
        \left(
            \sum_{i=1}^n
            \cc{w}^i a_i - w^i a_i^\dag
        \right)^r
        \psi
    \end{equation}
    converges unconditionally in the Hilbert space topology.
\end{corollary}

Of course we can also rewrite this in terms of the position and
momentum operators
\begin{equation}
    \label{eq:Position}
    Q_k = \frac{1}{2}(a_k + a_k^\dag) 
    = \frac{1}{2}\pi(z^k + \cc{z}^k)
\end{equation}
and
\begin{equation}
    \label{eq:Momentum}
    P_k = \frac{1}{2\I}(a_k -a_k^\dag) 
    = \frac{1}{2\I} \pi(z^k -\cc{z}^k),
\end{equation}
defined as unbounded symmetric operators on $\mathfrak{D}$. Then
\eqref{eq:UwcStarExp} shows that for $\psi \in \mathfrak{D}$ we have
the unconditionally convergent series
\begin{equation}
    \label{eq:ExpQ}
    \E^{\frac{\I\vec{p}\cdot \vec{Q}}{\hbar}}\psi
    =
    \sum_{r=0}^\infty \frac{1}{r!}
    \left(\frac{\I \vec{p} \cdot \vec{Q}}{\hbar}\right)^r \psi
\end{equation}
and
\begin{equation}
    \label{eq:ExpP}
    \E^{\frac{\I\vec{q}\cdot \vec{P}}{\hbar}}\psi
    =
    \sum_{r=0}^\infty \frac{1}{r!}
    \left(\frac{\I \vec{q} \cdot \vec{P}}{\hbar}\right)^r \psi
\end{equation}
in the Hilbert space topology, for all $\vec{q}, \vec{p} \in
\mathbb{R}^n$ substituting $w$ suitably.

Of course, this also follows by `Hilbert space techniques' from the
strong continuity of the representation $U$ and
Theorem~\ref{theorem:AnalyticVectors}. Note however, that the above
argument using the convergence of the star exponential is independent.

%
%

\begin{footnotesize}
%

\end{footnotesize}

\end{document}